\documentclass{article}
\usepackage[greek,american]{babel}
\usepackage{amssymb}
\usepackage{amsfonts}
\begin{document}
\newtheorem{remark}{Remark}[section]
\newtheorem{pigsfly1}{Definition}[section]
\newtheorem{rigor1}{Theorem}[section]
\newtheorem{rigor2}{Proposition}[section]
\newtheorem{smallTheorem}{Lemma}[section]
\newtheorem{help1}{Example}[section]
\title{Global Existence and Compact Attractors for the  Discrete Nonlinear Schr\"odinger  equation}
\author{Nikos I. Karachalios and Athanasios N. Yannacopoulos\\ Department of Statistics and Actuarial Science\\ University of the Aegean \\Karlovassi, 83200, Samos, Greece
}

\maketitle
\begin{abstract}
We study the asymptotic behavior of solutions of discrete nonlinear Schr\"odinger-type (DNLS) equations. For a conservative system, we consider the global in time solvability and the question of existence of standing wave solutions. Similarities and differences with the continuous counterpart (NLS-partial differential equation) are pointed out.  For a dissipative system we prove existence of a global attractor and its stability under finite dimensional approximations. Similar questions are treated in a weighted phase space. Finally, we propose possible extensions for various types of DNLS equations.
\end{abstract}
\section{Introduction}
The discrete nonlinear Schr\"odinger equation (DNLS), is a very popular model with a great variety of applications, ranging from physics to biology. To name but a few, the DNLS has been successfully applied to the modelling of localized pulse propagation in optical fibers and wave guides, to the study of energy relaxation in solids, the behaviour of amorphous material, to the modelling of self-trapping of vibrational energy in proteins or studies related to the denaturation of the DNA double strand (for an account of possible applications see e.g. the recent review \cite{Kevrekidis}). The DNLS serves both as a model in its own right (modelling cases where the nature of the problem is inherently discrete), or as an approximation of the continuous nonlinear Schr\"odinger through a numerical scheme.

The study of the dynamics of the discrete nonlinear Schr\"odinger equation, has been an active theme of research in the past decade. The major part of the activity, has focused on the study of localized excitations and nonlinear waves (solitons and breathers), in the conservative DNLS. This activity has started with the ingenious construction introduced in \cite{MackayAubry}, of localized in space time periodic or time-quasiperiodic solutions of dicrete lattice systems (which include the DNLS as special case), using the so called anti-integrable limit or anti-continuous limit. Such solutions are essentially found for the weak coupling limit between the lattice
sites. Other constructions are available using dynamical systems methods, such as the reduction to a symplectic map and the homoclinic trajectory method. The stability of such solutions has been studied using both analytic and numerical techniques and they have been shown in many cases to be robust under perturbations. They can be appplied to the study of diverse physical and biological phenomena such as lattice dynamics of solids, selective bond excitations, localization of electromagnetic waves in photonic crystals with nonlinear response etc. For detailed reviews on this theme we refer to  \cite{Aubry},  \cite{FlachWillis}, \cite{HennigTsironis}. There have also been recent interesting works on more complicated solutions arising in higher dimensions than one, such as discrete vortex solutions (see eg. \cite{Kevrekidis1} or \cite{Kevrekidis2}).

In this work we mainly consider the  asymptotic behavior of solutions for the following complex lattice dynamical system (LDS)
\begin{eqnarray}
\label{lat1}
i\dot{u}_n+\frac{1}{\epsilon}(u_{n-1}-2u_n+u_{n+1}) &+&i\delta u_n+F(u_n)=g_n,\;\;n\in\mathbb{Z},\;\;\epsilon >0,\\
\label{lat2}
u_n(0)&=&u_{n,0},\;\;n\in\mathbb{Z}
\end{eqnarray}
The real parameter $\delta >0$ introduces weak damping, while the case $\delta =0$ and $g_n=0,$ for all $n\in\mathbb{Z}$, corresponds to a conservative
(NLS) lattice system.

Using techniques from the theory of infinite dimensional systems,
we are able to show the global in time existence of solutions to
both the conservative and the dissipative DNLS equation. This
ensures the well posedeness of the model. This result is
interesting as it provides a rigorous proof of the fact that blow
up is not expected to occur in the discrete models as it may happen
in the continuous model, a fact that has been observed numerically
by a number of authors \cite{bang}. Futhermore, for a finite lattice (assuming Dirichlet boundary conditions), and by using a discrete
version of the mountain pass Theorem,  we prove the existence of nontrivial
standing waves, in the conservative case. For the infinite dimensional case, we show non-existence of nontrivial standing waves for certain values of energy, by using a fixed point argument. For the
dissipative case we may further prove a more strong result; the
existence of a global attractor, that attracts all bounded sets of
the infinite dimensional phase space. It can further be shown that
the dynamics on this global attractor may be approximated by the
dynamics of finite dimensional systems. Results of this kind, for
lattice dynamical systems of first and second order are presented
in the recent works \cite{Bates,SZ1,SZ2}.

Many of the physically interesting solutions of the DNLS equation present strong spatial localization properties. This is true for instance for soliton solutions or breathers. To cover such situations, we study the DNLS equation in weighted spaces, with properly chosen weight functions. For the dissipative case we may prove global existence of solutions in weighted spaces as well as the existence of a global attractor. Finally, we comment on possible generalizations of the above results in higher spatial dimensions and for more general models.
%%%%%%%%%%%%%%%%%%%%%%%%%%%%%%%%%%%%%%%%%%%%%%%%%%%%%%%%555%%%%%%%%%%%%5
%%%%%%%%%%%%INTRO%%%%%%%%%%%%%%%%%%%%%%%%%%%%%%%%%%%%%%%%%%%%%%%%%%%%%%%
%%%%%%%%%%%%%%%%%%%%%%%%%%%%%%%%%%%%%%%%%%%%%%%%%%%%%%%%%%%%%%%%%%%%%%%%
%%%%%%%%%%%%%%%%%%%%%%%%%%%%%%%%%%%%%%%%%%%%%%%%%%%%%%%%%%%%%%%%%%%%%%%%
%%%%%%%%%%%%%%%%%%%%%%%%%%%%%%%%%%%%%%%%%%%%%%%%%%%%%%%%%%%%%%%%%%%%%%%%
\section{Formulation of the problem}
Problem (\ref{lat1})-(\ref{lat2}) is an infinite system of ordinary differential
equations (infinite dimensional dynamical system). We  will  formulate the problem  in an appropriate functional setting.
In what follows, a
complex Hilbert space $\mathrm{X}$ with the sesquilinear form
$B_{\mathrm{X}}(\cdot ,\cdot)$ and the norm
$||\cdot||_{\mathrm{X}}$, will be considered as a {\em real}
Hilbert space with the scalar product $(\cdot
,\cdot)_{\mathrm{X}}=\mathrm{Re}\,B_{\mathrm{X}}(\cdot
,\cdot)$.Let
$\mathbf{T}:D(\mathbf{T})\subseteq X\rightarrow X$  a $\mathbb{C}$-linear, self-adjoint\ $\leq 0$ operator with dense domain $D(\mathbf{T})$ on the Hilbert space $X$, equipped with the scalar product $(\cdot ,\cdot)_{X}$. The space
$X_{\mathbf{T}}$ is the completion of $D(\mathbf{T})$ in the norm $||u||_{\mathbf{T}}^2=||u||^2_X-(\mathbf{T}u,u)_X$, for $u\in X_{\mathbf{T}}$, and we denote
by $X_{\mathbf{T}}^*$ its dual and by $\mathbf{T}^*$ the extension of $\mathbf{T}$ to the dual of $D(T)$, denoted by $D(T)^*$. A function $F:X\rightarrow X$  is Lipschitz
continuous on bounded subsets of $X$, if for all $M>0$ there exists a constant $L(M)$ such that $||F(u)-F(v)||_X\leq L(R)||u-v||_X$, for all $u,v\in B_M$ the closed ball of $X$, of center $0$ and of radius $R$.

For the treatment of (\ref{lat1})-(\ref{lat2}), we shall use complexifications  of the usual real sequence spaces denoted by
\begin{eqnarray}
\label{lp}
{\ell}^p:=\left\{u=(u_n)_{n\in\mathbb{Z}}\in\mathbb{C}\;\;:\;\;
||u||_{\ell^p}:=\left(\sum_{n\in\mathbb{Z}}|u_n|^p\right)^{\frac{1}{p}}<\infty\right\}.
\end{eqnarray}
Let us recall that between $\ell^p$, spaces the following elementary embedding relation  \cite[pg. 145]{HiLa} holds,
\begin{eqnarray}
\label{lp1}
\ell^q\subset\ell^p,\;\;\;\; ||u||_{\ell^p}\leq ||u||_{\ell^q}\,\;\; 1\leq q\leq p\leq\infty.
\end{eqnarray}
We will see in the sequel that relation (\ref{lp1}) has important effects concerning the global in time solvability of DNLS equation, especially in the conservative case.  For $p=2$ we get the usual Hilbert space of square-summable (complex) sequences endowed with the real scalar product
\begin{eqnarray}
\label{lp2}
(u,v)_{\ell^2}=\mathrm{Re}\sum_{{n\in\mathbb{Z}}}u_n\overline{v_n},\;\;u,\,v\in\ell^2.
\end{eqnarray}
We describe next the discrete differential operators employed. For any $u,v\in\ell^2$
we consider  the linear operators $A,B,B^*:\ell^2\rightarrow\ell^2$,
\begin{eqnarray}
\label{diffop}
(Bu)_{n\in\mathbb{Z}}&=&u_{n+1}-u_{n},\;\;\;\;\;\;
(B^*u)_{n\in\mathbb{Z}}=u_{n-1}-u_{n},\\
(Au)_{n\in\mathbb{Z}}&=&(u_{n-1}-2u_n+u_{n+1}).
\end{eqnarray}
Moreover it holds that
\begin{eqnarray}
\label{diffop2}
(Bu,v)_{\ell^2}&=&(u,B^*v)_{\ell^2},\;\;\;\\
\label{diffop3}
(Au,v)_{\ell^2}&=&-(Bu,Bv)_{\ell^2},\;\;u,v\in\ell^2.
\end{eqnarray}
The functional set up will also include (see \cite{SZ1}) the following bilinear form and induced norm
\begin{eqnarray}
\label{lp3}
(u,v)_{\ell^2_1}&:=&(Bu,Bv)_{\ell^2}+ (u,v)_{\ell^2},\\
\label{lp4}
||u||_{\ell^2_1}^2&:=&||Bu||_{\ell^2}^2+||u||_{\ell^2}^2.
\end{eqnarray}
We denote by $\ell^2_1$ the Hilbert space with  scalar product (\ref{lp3}) and norm  (\ref{lp4}). The usual norm of $\ell^2$ and (\ref{lp4}) are equivalent, since
\begin{eqnarray}
\label{lp5}
||u||^2_{\ell^2}\leq ||u||^2_{\ell^2_1}\leq c||u||^2_{\ell^2},\;\;c>0.
\end{eqnarray}

For local existence, we examine the following examples of nonlinearities $F:\mathbb{C}\rightarrow\mathbb{C}$: \vspace{.2cm}
\newline
$\mathrm{(N_1)}$\ \ $F(0)=0$ and there exist constants $c>0$, $\beta\geq 0$ such that
$|F(z_1)-F(z_2)|\leq c(|z_1|^{\beta}+|z_2|^{\beta})|z_1-z_2|$, or alternatively
\newline
$\mathrm{(N_2)}$\ \ $F(z)=f(|z|^2)z$ where $f,\,f':\mathbb{R}\rightarrow\mathbb{R}$, are continuous.
\newline

Both cases include the power-law nonlinearity
$F(z)=|z|^{2\sigma}z$ (for $\mathrm{(N_2)}$ we must have  $1<\sigma <\infty$). Since we
intend to apply the  general theory of abstract Schr\"{o}dinger
equations \cite{cazh,cazS}, the following Lemma will be useful.
\begin{smallTheorem}
\label{LipN}
Let $F:\mathbb{C}\rightarrow\mathbb{C}$ satisfy $\mathrm{(N_1)}$ or $\mathrm{(N_2)}$. Then the function $F$ defines an operator (still denoted by $F$) $$F:\ell^2\rightarrow\ell^2,\;\;(F(u))_{n\in\mathbb{Z}}:=F(u_n),$$ which is Lipschitz continuous on bounded sets of $\ell^2$.
\end{smallTheorem}
{\bf Proof:}\ \ We consider first the case $\mathrm{(N_1)}$.  Let $u\in B_R$ a closed ball in $\ell^2$ of center $0$ and radius $R$.  We have from (\ref{lp}) that
\begin{eqnarray}
\label{ufp1}
||\mathrm{F}(u)||^2_{\ell^2}&\leq& c^2\sum_{n\in\mathbb{Z}}|u_n|^{2\beta+2}
=c^2||u||_{\ell^{2\beta +2}}^{2\beta +2}\leq c^2||u||_{\ell^2}^{2\beta +2}
\end{eqnarray}
hence $\mathrm{F}:\ell^2\rightarrow\ell^2$,  is  bounded on bounded
sets of $\ell^2$.

For $u,v\in B_R$, we observe that
\begin{eqnarray}
\label{ufp2}
||\mathrm{F}(u)-F(v)||^2_{\ell^2}&\leq& c^2\sum_{n\in\mathbb{Z}}(|u_n|^{\beta}+|v_n|^{\beta})^2|u_n-v_n|^2\\
&\leq& c^2\sup_{n\in\mathbb{Z}}\left[(|u_n|^{\beta}+|v_n|^{\beta})^2\right]
\sum_{n\in\mathbb{Z}}|u_n-v_n|^2\leq c^24R^{2\beta}||u-v||_{\ell^2}^2,
\end{eqnarray}
justifying that the map $F:\ell^2\rightarrow\ell^2$ is
Lipschitz continuous on bounded sets of $\ell^2$ with Lipschitz constant $L(R)=c2R^{\beta}$.\newline

For the case $\mathrm{(N_2)}$ we have
\begin{eqnarray}
\label{prop1}
||\mathrm{F}(u)||^2_{\ell^2}=\sum_{n\in\mathbb{Z}}|f(|u_n|^2)|^2|u_n|^2.
\end{eqnarray}
Since $f:\mathbb{R}\rightarrow\mathbb{R}$ is continuous, there exists a monotone increasing $\mathrm{C^1}$-function
$g:\mathbb{R}\rightarrow\mathbb{R}$ such that
\begin{eqnarray}
\label{prop2}
|f(\rho)|\leq g(|\rho|),\;\;\mbox{for all}\;\;\rho\in\mathbb{R}.
\end{eqnarray}
(see e.g \cite[p.g 796]{zei85}). Note also that $|u_n|^2\leq ||u_n||^2_{\ell^2}$ for all $n\in\mathbb{Z}$. Now (\ref{prop1})
and (\ref{prop2}) imply the inequality
\begin{eqnarray}
\label{prop3}
||\mathrm{F}(u)||_{\ell^2}^2\leq \sum_{n\in\mathbb{Z}}g(|u_n|^2)^2|u_n|^2
&\leq& \sum_{n\in\mathbb{Z}}g(||u||_{\ell^2}^2)^2|u_n|^2\nonumber\\
&\leq&\{\mathrm{max}_{\rho\in [0,R^2]}g(\rho)\}^2\sum_{n\in\mathbb{Z}}|u_n|^2\nonumber\\
&\leq&K(R)^2||u||^2_{\ell^2}.
\end{eqnarray}
with $K(R)=g(R^2)$.

To check  the Lipschitz property, we may argue as for the proof of (\ref{prop3}).  For some $\theta\in (0,1)$ and an appropriate $\mathrm{C^1}$-function
$g_1:\mathbb{R}\rightarrow\mathbb{R}$, we get the inequality
\begin{eqnarray}
\label{prop4}
||\mathrm{F}(u)-\mathrm{F}(v)||^2_{\ell^2}
&\leq&2\sum_{n\in\mathbb{Z}}|f(|u_n|^2)|^2|u_n-v_n|^2
+2\sum_{n\in\mathbb{Z}}|f(|u_n|^2)-f(|v_n|^2)|^2|v_n|^2,\\
\label{prop5}
\sum_{n\in\mathbb{Z}}|f(|u_n|^2)-f(|v_n|^2)|^2|v_n|^2&=&
\sum_{n\in\mathbb{Z}}|f'(\theta |u_n|^2+(1-\theta)|v_n|^2)|^2(|u_n|+|v_n|)^2
(|u_n|-|v_n|)^2|v_n|^2\nonumber\\
&\leq&\{\mathrm{max}_{\rho\in [0,2R^2]}g_1(\rho)\}^2\sup_{n\in\mathbb{Z}}\left[(|u_n|+|v_n|)^2|v_n|^2\right]
\sum_{n\in\mathbb{Z}}|u_n-v_n|^2.
\end{eqnarray}
A combination of (\ref{prop4}) and (\ref{prop5}), implies that there exists a constant
$L(R)=(2g^2(R^2)+8R^4g_1(2R^2))^{1/2}$ such that
\begin{eqnarray}
\label{prop6}
||\mathrm{F}(u)-\mathrm{F}(v)||_{\ell^2}\leq L(R)||u-v||_{\ell^2}.
\end{eqnarray}
The Lemma is proved. \ \ $\diamond$\newline

We observe by (\ref{diffop3}), that the operator $A$ satisfies the relations
\begin{eqnarray}
\label{lp6}
(Au,u)_{\ell^2}&=&-||Bu||_{\ell^2}^2\leq 0,\\
\label{lp7}
(Au,v)_{\ell^2}&=&(u,Av)_{\ell^2},
\end{eqnarray}
therefore defines a self-adjoint operator on $D(A)=X=\ell^2$ and $A\leq 0$. Note that the graph norm
$$||u||_{D(A)}=||Au||^2_{\ell^2}+||u||^2_{\ell^2},$$ is also an equivalent norm with the $\ell^2$-norm since
$$||u||_{\ell^2}^2\leq\sum_{n\in\mathbb{Z}}|u_{n+1}-2u_n+u_{n-1}|^2+\sum_{n\in\mathbb{Z}}|u_n|^2\leq c||u||_{\ell^2}^2.$$
In our case, as it is indicated by
(\ref{lp3})-(\ref{lp6})-(\ref{lp7}), we may choose $X_A=\ell^2_1$
equipped with the norm $||u||_A^2=||u||_X^2-(Au,u)_X\equiv
||u||_{\ell^2_1}$, for $u\in\ell^2$. Moreover,
$D(A)=X=\ell^2=D(A)^*$. Obviously $A^*=A$ and the operator
$iA:\ell^2\rightarrow \ell^2$ defined by $(iA)u=iAu$ for $u\in
\ell^2$, is $\mathbb{C}$-linear and skew-adjoint and $iA$
generates a group $(\mathcal{T}(t))_{t\in\mathbb{R}}$, of
isometries on $\ell^2$. It is easy to check that the same properties hold for the operator $\epsilon^{-1} A$ with $\epsilon>0$ (i.e., we may consider the equivalent norm $\epsilon^{-1}||Bu||^2_{\ell^2}+||u||_{\ell^2}$). 

Thus, for fixed $T>0$ and $(u_{n,0})_{n\in\mathbb{Z}}:=u_0\in\ell^2$, a function $u\in\mathrm{C}([0,T],\ell^2)$ is a solution of (\ref{lat1})-(\ref{lat2}) if and only if
\begin{eqnarray}
\label{milds}
u(t)=\mathcal{T}(t)u_0+i\int_{0}^{t}\mathcal{T}(t-s)F_1(u(s))ds,\;\;
(F_1(u))_{n\in\mathbb{Z}}:=F_1(u_n)=i\delta u_n+F(u_n)-g_n.
\end{eqnarray}
Our local existence result can be stated as follows:
\begin{rigor1}
\label{locex}
Let assumptions $(\mathrm{N_1})$ or $(\mathrm{N_2})$ be satisfied, assume that $g:=(g_n)_{n\in\mathbb{Z}}\in\ell^2$. Then there exists a function $T^*:\ell^2\rightarrow (0,\infty]$ with the following properties:\vspace{.2cm}\\
(a)\  For all $u_0\in\ell^2$, there exists  $u\in\mathrm{C}([0,T^*(u_0)),\ell^2)$ such that
for all $0<T<T^*(u_0)$, $u$ is the unique solution of (\ref{lat1})-(\ref{lat2}) in $\mathrm{C}([0,T],\ell^2)$ (well posedeness).\vspace{.2cm}\\
(b)\  For all $t\in [0,T^*(u_0))$,
\begin{eqnarray}
\label{maxT}
T^*(u_0)-t\geq \frac{1}{2(L_1(R)+1)}:=T_R,\;\;R=||g||_{\ell^2}+2||u(t)||_{\ell^2},\;\;L_1(R)=\delta+L(R),
\end{eqnarray}
where $L(R)$ is the Lipschitz constant for the map $F:\ell^2\rightarrow\ell^2$. Moreover the following alternative holds: (i) $T^{*}(u_0)=\infty$, or (ii) $T^*(u_0)<\infty$ and $\lim_{t\uparrow T^*(u_0)}||u(t)||_{\ell^2}=\infty$ (maximality).\vspace{.2cm}\\
(c)\ $T^*:\ell^2\rightarrow (0,\infty]$ is lower semicontinuous.
In addition, if $\{u_{n0}\}_{n\in\mathbb{N}}$ is a sequence in
$\ell^2$ such that $u_{n0}\rightarrow u_0$ and if $T<T^*(u_0)$,
then $S(t)u_{0n}\rightarrow S(t)u_0$ in $\mathrm{C}([0,T],\ell^2)$, where $S(t)u_0=u(t)$, $t\in [0,T^*(u_0))$, denotes the solution operator
(continuous dependence on initial data).
\end{rigor1}
{\bf Proof:}\ Since the result is an application of the results in
\cite[p.g. 56-59]{cazh}, we present only an outline of the proof,
for the sake of completeness. It follows from Lemma \ref {LipN}
that $F_1:\ell^2\rightarrow\ell^2$ is Lipschitz continuous on
bounded sets of $\ell^2$ with Lipschitz constant
$L_1(R)=\delta+L(R)$. For all $u_0\in\ell^2$ with
$||u_0||_{\ell^2}\leq M$ we set
$R=2M+||F_1(0)||_{\ell^2}=2M+||g||_{\ell^2}$. Following closely
the lines of \cite[Lemma 4.3.2-Proposition 4.3.3]{cazh} and using
(\ref{milds}) we may show that the map
$$\Phi_u(t):=\mathcal{T}(t)u_0+i\int_{0}^{t}\mathcal{T}(t-s)F_1(u(s))ds,$$
is a contraction on the complete metric space $(X_R,d)$ where
$$X_R:=\left\{u\in\mathrm{C}([0,T_R],\ell^2):||u(t)||_{\ell^2}\leq
R,\;\forall t\in [0,T_R]\right\},\;\;
T_R=\frac{1}{2L_1(R)+1},\;\;d=\max_{t\in
[0,T_R]}||u(t)-v(t)||_{\ell^2}.$$ For $u_0\in\ell^2$ we define
$T^*(u_0)=\sup\left\{T>0:\exists\;\; S(t)u_0=u(t)\in
\mathrm{C}([0,T],\ell^2)\;\;\mbox{ solution of}\;\;
(\ref{milds})\right\}$. Inequality (\ref{maxT}) follows by the
contradiction argument of \cite[Theorem 4.3.4]{cazh}: We assume
instead of (\ref{maxT}), that for $R=||u(t)||_{\ell^2}$,
$T^*(u_0)-t<T_R$ and we consider the solution
$v\in\mathrm{C}([0,T_R],\ell^2)$ of
$$v(s)= \mathcal{T}(s)u(t)+i\int_{0}^{s}\mathcal{T}(s-\sigma)F_1(v(\sigma))d\sigma,\;\;s\in [0,T_R].$$ Then we can define a solution $w\in\mathrm{C}([0,t+T_R],\ell^2)$ of (\ref{milds}) at $T=t+T_R$, defined by
\begin{eqnarray*}
w(s)=\left\{
\begin{array}{ll}
u(s),\;\;&s\in [0,t], \\
v(s-t),\;\;&s\in [t,t+T_R].
\end{array}
\right.
\end{eqnarray*}
which is in contradiction with the assumption $T^*(u_0)<t+T_R$. To
prove (c), as in \cite[Proposition 4.3.7]{cazh}, we set
$R=2\sup_{t\in[0,T]}||u(t)||_{\ell^2}$ and define
$\tau_n=\sup\left\{t\in
[0,T^*(u_{0n})):||S(s)u_{0n}||_{\ell^2}\leq 2R,\;\forall s\in
[0,t]\right\}$. Since $u_{0n}\rightarrow u_0$ in $\ell^2$, it
follows that for sufficiently large $n$, $||u_{0n}||_{\ell^2}<R$.
Therefore $\tau_n>T_R>0$. Moreover, for all $t\leq T$ and
$t\leq\tau_n$, we get from (\ref{milds}) and Gronwall's Lemma the
inequality
\begin{eqnarray*}
||S(t)u_{0n}-S(t)u_{0}||_{\ell^2}\leq
e^{L_1(2R)T}||u_{0n}-u_{0}||_{\ell^2}.
\end{eqnarray*}
Letting $n\rightarrow\infty$, we get that
$S(t)u_{0n}\rightarrow S(t)u_{0}$ in $\ell^2$ and $\tau_n>T$, which
implies that $T^*(u_{0n})>T$. \ $\diamond$
\section{Global existence for the conservative case}
In this section we discuss the global in time solvability of (\ref{lat1})-(\ref{lat2}) for the case $\delta=0$ and $g_n=0$, for all $n\in\mathbb{Z}$ (conservative case). We consider as a model problem, the following DNLS equation with arbitrary power nonlinearity
\begin{eqnarray}
\label{cons1}
i\dot{u}_n+\frac{1}{\epsilon}(u_{n-1}&-&2u_n+u_{n+1}) +|u_n|^{2\sigma}u_n=0,\;\;n\in\mathbb{Z},\;\;0\leq\sigma <\infty,\;\;\epsilon>0,\\
\label{cons2}
u_n(0)&=&u_{n,0},\;\;n\in\mathbb{Z}
\end{eqnarray}
We also consider the
 (NLS) partial differential equation
\begin{eqnarray}
\label{pde}
i\partial_t u + u_{xx}&+&|u|^{2\sigma}u=0,\;\; x\in\mathbb{R},\;\;t>0,\\
u(x,0)&=&u_0(x).\nonumber
\end{eqnarray}

We will show that there is a vast difference concerning the global
solvability of the discrete (NLS) equation (\ref{cons1}) with $\epsilon >0$ and its
continuum limit as $\epsilon\rightarrow 0$ given by (\ref{pde}).
Concerning (\ref{cons1}), we have the following Theorem.
\begin{rigor1}
\label{global}
Let $u_0\in\ell^2$. For any $\epsilon >0$, and $0\leq\sigma <\infty$, there exists a unique  solution of (\ref{cons1}) such that $u(t)\in \mathrm{C}^1([0,\infty),\ell^2)$.
\end{rigor1}
{\bf Proof:}\ \  The result is a consequence of discrete conservation laws satisfied by
the solutions of (\ref{cons1})-(\ref{cons2}). We take the scalar product of (\ref{cons1}) with $iu$. By using (\ref{lp6}) we obtain
\begin{eqnarray}
\label{charge}
\frac{d}{dt}||u(t)||^2_{\ell^2}=0,\;\;\mbox{or}\;\;||u(t)||^2_{\ell^2}=||u_0||^2_{\ell^2},\;\;\mbox{for every}\;\;t\in [0,T^*(u_0)),
\end{eqnarray}
i.e. $||u(t)||_{\ell^2}$ is uniformly bounded on the maximal
interval of existence. It follows then by Theorem \ref{locex}
(b)-(c) that $T^*(u_0)=\infty$ and
$\mathrm{sup}\left\{||u(t)||_{\ell^2},\,t\in[0,\infty)\right\}<\infty$.\
\ $\diamond$
\begin{remark}
\label{ENERGY}
{\em Theorem \ref{global} can be also established by  a discrete version of the conservation of energy, satisfied by solutions of (\ref{cons1}). This approach elucidates the role of the nonlinearity exponent, in the global in time solvability of DNLS: It can be easily checked that 
\begin{eqnarray}
\label{eneg2}
\mathrm{E}(u(t))=\mathrm{E}(u_0),\;\;\mathrm{E}(u(t)):=\frac{1}{\epsilon}||u(t)||^2_{\ell^2_1}-\frac{1}{\sigma +1}||u(t)|_{\ell^{2\sigma +2}}^{2\sigma +2},\;\;t\in [0,T^*(u_0)).
\end{eqnarray}
Then, by using   (\ref{lp1}) and (\ref{eneg2}) we may derive the estimate
\begin{eqnarray}
\label{lastest}
||u(t)||_{\ell^{2}_1}^2
\leq  ||u_0||_{\ell^2_1}^2+ \frac{2\epsilon}{\sigma +1}||u_0||^{2\sigma +2}_{\ell^2}.
\end{eqnarray}
As a consequence of (\ref{lastest}) we obtain once again that $T^*(u_0)=\infty$ and $\mathrm{sup}\left\{||u(t)||_{\ell^2_1},\,t\in [0,\infty)\right\}<\infty$.

Clearly, Theorem \ref{global}, covers the case of DNLS with a
nonlinearity satisfying $\mathrm{(N_2)}$. The case of a globally
Lipschitz continuous function $F:[0,\infty)\rightarrow\mathbb{R}$
with $F(0)=0$ is also included. Following \cite[Section
7.2]{cazh}, \cite[p.g 53]{cazS}, such a function can be extended
to the complex function
$F:\mathbb{C}\setminus\{0\}\rightarrow\mathbb{C}$ by setting
$F(z)=\frac{z}{|z|}F(|z|)$. It can be shown as in Lemma \ref{LipN}
that $F:\ell^2\rightarrow\ell^2$ is locally Lipschitz and
$(F(u),iu)_{\ell^2}=0$ which implies (\ref{charge}). Setting
$F_*(z)=\int_{0}^{|z|}F(s)ds$ it can be shown as in Section 4,
Lemma \ref{derivative} that the functional
$(V(u))_n:=-\sum_{n\in\mathbb{Z}}F_*(u_n)$ is a
$\mathrm{C}^1$-functional on $\ell^2$ and $V'(u)=-F(u)$. In this
case, we observe conservation of energy defined by
$\mathrm{E}(u(t)):=\frac{1}{\epsilon}||u(t)||^2_{\ell^2_1}+V(u)$.\
\ $\bullet$}
\end{remark}

For a comparison  between the discrete  (NLS) equation (\ref{cons1}) and its continuous counterpart (\ref{pde}), we refer to the main results concerning (\ref{pde})
(see \cite{cazh,cazS,Yvan}).

For $u_0\in\mathrm{H^1}(\mathbb{R})$ and $0\leq\sigma <\infty$ there exists a unique maximal solution of (\ref{pde}), $u(t)\in \mathrm{C}([0,T_{max}),\mathrm{H^1}(\mathbb{R}))\cap\mathrm{C}^1([0,T_{max}),\mathrm{L^2}(\mathbb{R}))$. In addition:
If $\sigma <2$ then $T_{\max}=\infty$ and $u$ is bounded in $\mathrm{H^1}(\mathbb{R})$.
Let $\sigma\geq 2$. Assume that $u_0\in \mathrm{H^1}(\mathbb{R})$ such that $\int_{\mathbb{R}}|x|^2|u_0|^2dx <\infty$ (initial data with finite variance) and $E(u_0)<\infty$. Then $T_{\max}<\infty$.
On the other hand if $||u_0||_{\mathrm{H^1}}$ is sufficiently small, $T_{max}=\infty$ and $u$ is bounded in $\mathrm{H^1}(\mathbb{R})$.

It follows from Theorem \ref{global}, that solutions of (\ref{cons1}) with $\epsilon >0$, exist globally in $\ell^2_1$-norm unconditionally with respect to the degree of the nonlinearity the size of the initial data and the sign of the initial energy $\mathrm{E}(u_0)$.  This is in agreement with the numerical observations in \cite{bang}: According to \cite{bang}, for a discrete system of the form (\ref{cons1}), numerical simulations provide evidence that the solution of the discrete equation still exists after localization, while
the solution of the continuum system blows-up in finite time. Note that the discrete version of initial data with finite variance reads as $\sum_{n\in\mathbb{Z}}n^2|u_{n,0}|^2<\infty$ and such data belong to $\ell^2$.
%%%%%%%%%%%%%
%%%%%%%%%%%%%
%%%%%%%%%%%%%%
\section{Existence of standing wave solutions for the conservative DNLS}
We conclude the discussion on the conservative DNLS,  with a discussion on the existence of standing wave solutions. The standing wave solution for fixed $\omega^2$ is given by the ansatz
\begin{eqnarray}
\label{sw}
u_n(t)= e^{i\omega^2t}\phi_n,\;\;n\in\mathbb{Z},\;\;\phi_n\;\;\mbox{independent of}\;\;t.
\end{eqnarray}
We do not restrict $\phi_n$ to be real.  We are focused on DNLS (\ref{cons1}). It can be easily seen, that any
standing wave solution of (\ref{cons1}), satisfies the following equation
\begin{eqnarray}
\label{swe}
-\frac{1}{\epsilon}(\phi_{n-1}-2\phi_n+\phi_{n+1})+\omega^2\phi_n=|\phi_n|^{2\sigma}\phi_n,\;\;n\in\mathbb{Z}.
\end{eqnarray}

We will study in this section two related problems. 
The first problem  we will study is given by equation (\ref{swe}), with $n$ taking finite values.
We consider the finite dimensional subspace of $\ell^2$
$$\mathbb{C}^{2m+1}:=\left\{\psi_n\in\ell^2\;\;:\;\;\psi_{-(m+1)}=\psi_{m+1}=0\right\}.$$ We consider the following nonlinear system in  $\mathbb{C}^{2m+1}$,
\begin{eqnarray}
\label{swef}
-\frac{1}{\epsilon}(\phi_{n-1}-2\phi_n&+&\phi_{n+1})+\omega^2\phi_n=|\phi_n|^{2\sigma}\phi_n,\;\; |n|\leq m,\\
\label{swefc}
\phi_{-(m+1)}&=&\phi_{m+1}=0.
\end{eqnarray}
This is a finite dimensional problem, and is related to the problem of existence of standing wave solutions for the DNLS equation in a finite lattice with Dirichlet boundary conditions. Although this problem is finite dimensional, it is still of interest as any numerical approximation of the continuous space NLS, or even of the infinite dimensional DNLS will necessarily lead to a finite dimensional problem of the above type. 

The set $\mathbb{C}^{2m+1}$ endowed with the discrete inner product and induced norm
\begin{eqnarray}
\label{discn}
(\phi,\psi)_2:=\mathrm{Re}\sum_{n=-m}^{n=m}\phi_n\overline{\psi_n},\;\;||\psi||_2:=\sum_{n=-m}^{n=m}|\psi_n|^2,\;\;\phi,\,\psi\in \mathbb{C}^{2m+1},
\end{eqnarray}
is a (finite dimensional) Hilbert space. We consider now the operators
$$(A_1\psi)_{|n|\leq m}:=\psi_{n-1}-2\psi_n+\psi_{n+1},\;\;(B_1\psi)_{|n|\leq m}=\psi_{n+1}-\psi_n.$$
It can be easily checked (see also \cite[pg. 117]{Akriv}) that
\begin{eqnarray}
\label{byparts}
(-A_1\psi,\psi)_2=\sum_{n=-m}^{n=m}|\psi_{n+1}-\psi_n|^2,\;\;(-A_1\phi,\psi)_2=(B_1\phi,B_1\psi)_2.
\end{eqnarray}
Hence, we may also consider the inner product and the corresponding norm in $\mathbb{C}^{2m+1}$ 
\begin{eqnarray}
\label{discSob}
(\phi,\psi)_{1,2}:=(B_1\phi,B_1\psi)_2+(\phi,\psi)_2,\;\;||\psi||_{1,2}:=\sum_{n=-m}^{n=m}(|\psi_{n+1}-\psi_n|^2+|\psi_n|^2).
\end{eqnarray}
The norm in (\ref{discSob}) is equivalent with (\ref{discn}), and the constants in  the equivalence inequality, can be chosen in order to {\em be independent of} $m$.

The second is the problem given by equation (\ref{swe}), which is an infinite dimensional system. This is related with the problem of existence
of standing wave solutions of the conservative DNLS equation in the infinite lattice. 

In what follows we will use a variational principle to  show existence of nontrivial standing wave solutions for (\ref{swef})-(\ref{swefc}) and a fixed point argument to show non-existence of non-trivial standing waves, for certain parameter values for the infinite dimensional problem (\ref{swe}).
%%%%%%%%%%%%%%%%%%%%%%%%%%%%%%%%%%%%%%%%%%%%%%%%%%%%%%%%%%%%%%
%%%%%%%%%%%%%%%%%%%%%%%%%%%%%%%%%%%%
%%%%%%%%%%Subsection-Finite%%%%%%%%%%%%%%%%%%%%%%%%%%%%%%%%%
%%%%%%%%%%%%%%%%%%%%%%%%%%%%%%%%%%%%%%%%%%%%%%%%%%%%%%%%%%%%
\subsection{Existence of non trivial standing wave solutions for the finite dimensional problem}
We will treat this problem using a variational approach. 
We start with some observations and results which are common to both problems. Unless the opposite is explicitly stated, everything that follows will hold for both the finite dimensional and the infinite-dimensional system.

Solutions of (\ref{swe}) are critical points of the functional
\begin{eqnarray}
\label{Enegfun}
\mathbf{E}(\phi )=\frac{1}{2\epsilon}\sum_{n\in\mathbb{Z}}|(B\phi)_n|^2+\frac{\omega^2}{2}\sum_{n\in\mathbb{Z}}|\phi_n|^2-\frac{1}{2\sigma +2}\sum_{n\in\mathbb{Z}}|\phi_n|^{2\sigma +2}.
\end{eqnarray}
To establish differentiability of the functional
$\mathbf{E}:\ell^2\rightarrow\mathbb{R}$, we shall use the
following discrete version of the dominated convergence Theorem,
provided by \cite{Bates2}.
\begin{smallTheorem}
\label{dc}
Let $\{\psi_{i,k}\}$ be a double sequence of summable functions  (i.e $\sum_{i\in\mathbb{Z}}|\psi_{i,k}|<\infty$) and $\lim_{k\rightarrow\infty}\psi_{i,k}=\psi_{i}$, for all
$i\in\mathbb{Z}$. If there exists a summable sequence $\{g_{i}\}$ such that $|\psi_{i,k}|\leq g_{i}$ for all $i,k$'s, we have that
$\lim_{k\rightarrow\infty}\sum_{i\in\mathbb{Z}}\psi_{i,k}=\sum_{i\in\mathbb{Z}}\psi_{i}$.
\end{smallTheorem}
We then have the following Lemma.
\begin{smallTheorem}
\label{derivative}
Let $(\phi_n)_{n\in\mathbb{Z}}=\phi\in\ell^{2\sigma+2}$ for some $0<\sigma <\infty$. Then the functional
$$\mathbf{V}(\phi)=\sum_{n\in\mathbb{Z}}|\phi_n|^{2\sigma +2},$$
is a $\mathrm{C}^{1}(\ell^{2\sigma +2},\mathbb{R})$ functional and
\begin{eqnarray}
\label{gatdev}
<\mathbf{V}'(\phi),\psi>=(2\sigma +2)\mathrm{Re}\sum_{n\in\mathbb{Z}}|\phi_n|^{2\sigma}\phi_n\overline{\psi_n},\;\;\psi=(\psi_n)_{n\in\mathbb{Z}}\in\ell^{2\sigma +2}.
\end{eqnarray}
\end{smallTheorem}
{\bf Proof:}\ \ We assume that $\phi,\,\psi\in\ell^{2\sigma +2}$.
Then for any $n\in\mathbb{Z}$, $0<s<1$, we get from  the mean
value Theorem
\begin{eqnarray}
\label{mv}
\frac{\mathbf{V}(\phi +s\psi)-\mathbf{V}(\psi)}{s}&=&\frac{1}{s}\mathrm{Re}\sum_{n\in\mathbb{Z}}\int_{0}^{1}\frac{d}{d\theta}|\phi_n +
\theta s\psi_n|^{2\sigma +2}d\theta\nonumber\\
&=&(2\sigma +2)\mathrm{Re}\sum_{n\in\mathbb{Z}}\int_{0}^{1}|\phi_n+s\theta\psi_n|^{2\sigma}
(\phi_n+s\theta\psi_n)\overline{\psi_n} d\theta.
\end{eqnarray}
For the rhs of (\ref{mv}), we have the estimate
\begin{eqnarray*}
(2\sigma +2)\sum_{n\in\mathbb{Z}}|\phi_n+\theta s\psi_n|^{2\sigma +1}|\psi_n|
\leq
(2\sigma +2)\sum_{n\in\mathbb{Z}}\left(|\phi_n|+|\psi_n|\right)^{2\sigma +1}|\psi_n|=\sum_{n\in\mathbb{Z}}z_n
\end{eqnarray*}
The sequence $(z_n)_{n\in\mathbb{Z}}$ is summable since
\begin{eqnarray*}
\sum_{n\in\mathbb{Z}}|z_n|\leq c\left(\sum_{n\in\mathbb{Z}}(|\phi_n|+|\psi_n|)^{2\sigma +2}\right)^{\frac{2\sigma +1}{2\sigma +2}}\left(\sum_{n\in\mathbb{Z}}|\psi_n|^{2\sigma +2}\right)^{\frac{1}{2\sigma +2}}.
\end{eqnarray*}
Letting $s\rightarrow 0$, the existence of the Gateaux derivative (\ref{gatdev}) of the functional $\mathbf{V}:\ell^{2\sigma +2}\rightarrow\mathbb{R}$,  follows from Lemma \ref{dc} (discrete dominated convergence).

We show next that the functional $\mathbf{V}':\ell^{2\sigma +2}\rightarrow\ell^{\frac{2\sigma +2}{2\sigma +1}}$ is
continuous. For $(\phi_n)_{n\in\mathbb{Z}}=\phi\in\ell^{2\sigma +2}$, we set $F_1(\phi_n)=|\phi_n|^{2\sigma}\phi_n$. We consider a sequence $\phi^m\in\ell^{2\sigma +2}$ such that $\phi^m\rightarrow \phi$ in $\ell^{2\sigma +2}$. We get then the inequality
\begin{eqnarray}
\label{hoin}
\left|\left<\mathbf{V}'(\phi_m)-\mathbf{V}'(\phi),\,\psi\right>\right|\leq c||F_1(\phi_m)-F_1(\phi)||_{\ell^{q}}||\psi||_{\ell^p},\;\;q=\frac{2\sigma +2}{2\sigma +1},\;\;p=2\sigma +2.
\end{eqnarray}
We denote by $(\phi_m)_n$ the $n$-th coordinate of the sequence $\phi_m\in\ell^2$. Since $F_1$ satisfies condition $\mathrm{(N_1)}$ with $\beta=2\sigma$, by setting $\Phi_n=|(\phi_m)_n|^{2\sigma}+|\phi_n|^{2\sigma}$, we get from
H\"{o}lder's inequality that
\begin{eqnarray*}
\sum_{n\in\mathbb{Z}}|F_1((\phi_m)_n)-F_1(\phi_n)|^{q}&\leq& c\sum_{n\in\mathbb{Z}}(\Phi_n)^q|(\phi_m)_n-\phi_n|^{q}\nonumber\\
&\leq& c\left(\sum_{n\in\mathbb{Z}}|(\phi_m)_n-\phi_n|^{2\sigma +2}\right)^{\frac{1}{2\sigma +1}}
\left(\sum_{n\in\mathbb{Z}}(\Phi_n)^{\frac{\sigma +1}{\sigma}}\right)^{\frac{2\sigma}{2\sigma +1}}\rightarrow 0,\;\mbox{as}\;m\rightarrow\infty.\;\;\diamond
\end{eqnarray*}

By  (\ref{lp6}),(\ref{lp7}) and Lemma (\ref{derivative}), it
follows that the functional $\mathbf{E}$ defined by
(\ref{Enegfun}) is $\mathrm{C}^1(\ell^2,\mathbb{R})$. For convenience, we recall \cite[Definition 4.1, pg. 130]{CJ}
(PS-condition) and \cite[Theorem 6.1, pg. 140]{CJ} or
\cite[Theorem 6.1, pg. 109]{struwe} (Mountain Pass Theorem of Ambrosseti-Rabinowitz \cite{Amb}). 
\begin{pigsfly1}
\label{condc} Let $X$ be a Banach space and $\mathbf{E}:X\rightarrow\mathbb{R}$ be $\mathrm{C}^1$. We say that
$\mathbf{E}$ satisfies condition $(PS)$ if, for any sequence $\{u_n\}\in X$ such that $|\mathbf{E}(u_n)|$ is bounded and $\mathbf{E}'(u_n)\rightarrow 0$ as $n\rightarrow\infty$,
there exists a convergent subsequence. If this condition is only satisfied in the region where $\mathbf{E}\geq\alpha >0$ (resp $\mathbf{E}\leq -\alpha <0$) for all $\alpha >0$, we say $\mathbf{E}$ satisfies condition $(PS^+)$ (resp. $(PS^-)$).
\end{pigsfly1}
\begin{rigor1}
\label{mpass}
Let $\mathbf{E}:X\rightarrow\mathbb{R}$ be $C^1$ and satisfy (a) $\mathbf{E}(0)=0$, (b) $\exists\rho >0$, $\alpha >0:\;||u||_X=\rho$ implies $\mathbf{E}(u)\geq\alpha$, (c) $\exists u_1\in X :\;||u_1||_X\geq\rho$ and
$\mathbf{E}(u_1)<\alpha$. Define $$\Gamma=\left\{\gamma\in \mathrm{C}^0([0,1],X):\;\gamma (0)=0,\;\;\gamma (1)=u_1\right\}.$$
Let $F_{\gamma}=\{\gamma(t)\in X:\;0\leq t\leq 1\}$ and $\mathcal{L}=\{F_\gamma :\;\gamma\in \Gamma\}$. If $\mathbf{E}$ satisfies condition $(PS)$, then $$\beta:=\inf_{F_{\gamma}\in \mathcal{L}}\sup\{\mathbf{E}(v):v\in F_{\gamma}\}\geq\alpha$$ is a critical point of the functional $\mathbf{E}$.
\end{rigor1}

We shall verify definition \ref{condc} and the assumptions of
Theorem \ref{mpass} to show that problem (\ref{swef})-(\ref{swefc}), has a
nontrivial solution.
\begin{rigor1}
\label{swee}
Consider the finite dimensional problem (\ref{swef})-(\ref{swefc}).
Let $0<\sigma<\infty$ and $\epsilon >0$. Then for any $\omega\neq 0$, there exists a non-trivial solution of
 (\ref{swef})-(\ref{swefc}).
\end{rigor1}
{\it Proof}\ \ ({\em Condition} $(PS)$) We observe  that the norm
\begin{eqnarray}
\label{ell2e}
||\phi||_{\ell^2_{\epsilon}}=\left(\frac{1}{\epsilon}||B\phi||_{\ell^2}^2+\omega^2|\phi|||_{\ell^2}^2\right)^{\frac{1}{2}},
\end{eqnarray}
is  equivalent with the norm of $\ell^2_1$ defined by (\ref{lp4}) since
\begin{eqnarray}
\label{lpno}
\min\left\{\frac{1}{\epsilon},\,\omega^2\right\}||\phi||^2_{\ell^2_1}\leq ||\phi||^2_{\ell^2_{\epsilon}}\leq
\max\left\{\frac{1}{\epsilon},\,\omega^2\right\}||\phi||^2_{\ell^2_1},
\end{eqnarray}
and is also equivalent  with the usual norm of $\ell^2$, as
it follows from (\ref{lp5}). Now we let a sequence $\phi_m$ of
$\ell^2$ be such that $|\mathbf{E}(\phi_m)|<M$ for some $M>0$ and
$\mathbf{E}'(\phi_m)\rightarrow 0$ as $m\rightarrow\infty$. By using
(\ref{Enegfun}) and Lemma \ref{derivative}, we observe that for $m$ sufficiently large
\begin{eqnarray}
\label{boundP.S}
M\geq \mathbf{E}(\phi_m)-\frac{1}{2\sigma +2}\left<\mathbf{E}'(\phi_m),\phi_m\right>=
\left(\frac{1}{2}-\frac{1}{2\sigma +2}\right)||\phi_m||^2_{\ell^2_{\epsilon}}.
\end{eqnarray}
Therefore the sequence $\phi_m$ is bounded. Thus, we may extract a subsequence, still denoted by $\phi_m$, such that
\begin{eqnarray}
\label{weakcon}
\phi_m\rightharpoonup \phi\;\;in\;\;\ell^2,\;\;\mbox{as}\;\;m\rightarrow\infty.
\end{eqnarray}

Focusing now on (\ref{swef})-(\ref{swefc}), we observe that since in the  finite dimensional space $\mathbb{C}^{2m+1}$ the weak convergence coincides with the strong,  $\phi_m$ converges strongly to $\phi$ in $\mathbb{C}^{2m+1}$. Therefore the functional $\mathbf{E}$ associated with the boundary value problem (\ref{swef})-(\ref{swefc}), satisfies condition $(PS)$.\newline
({\em Mountain Pass assumptions})\ \ We now check the conditions for the validity of the Mountain Pass Theorem for the functional $\mathbf{E}$. We use the same notations for the equivalent norms in the finite dimensional space  $\mathbb{C}^{2m+1}$, with those of the norms of the infinite dimensional spaces, having in mind that we deal with finite dimensional sums.

For every $\phi\in \mathbb{C}^{2m+1}$ it holds that 
\begin{eqnarray}
\mid\mid \phi\mid\mid_{\ell^{2\sigma +2}} \le \kappa_1\mid\mid \phi \mid\mid_{\ell^{2}_{\epsilon}}, \;\;\kappa_1=\frac{1}{\min\{1/\epsilon,\omega\}}.
\nonumber
\end{eqnarray}
The following inequality holds
\begin{eqnarray}
\mathbf{E}(\phi) &=& \frac{1}{2}\mid\mid \phi \mid\mid_{\ell^{2}_{\epsilon}}^2- \frac{1}{2\sigma+2} \mid\mid \phi \mid\mid_{\ell^{2\sigma +2}}^{2\sigma +2}\nonumber\\
&\geq&\frac{1}{2}\mid\mid
\phi \mid\mid_{\ell^{2}_{\epsilon}}^2-\frac{\kappa_1^{2\sigma+2}}{2\sigma+2}||\phi||_{\ell^{2}_{\epsilon}}^{2\sigma +2}.\nonumber
\end{eqnarray}
Choosing $\phi \in \mathbb{C}^{2m+1}$ such that $\mid\mid \phi \mid\mid_{\ell^2_{\epsilon}} =r$, we observe that if  $0<r < \left(\frac{\sigma +1}{\kappa_1^{2\sigma +2}}\right)^{\frac{1}{2\sigma}}$, then
\begin{eqnarray}
\mathbf{E}(\phi) \geq \alpha>0,\;\;\alpha=r^2\left(\frac{1}{2}-\frac{\kappa_1^{2\sigma+2}}{2\sigma +2}r^{2\sigma}\right).
\nonumber
\end{eqnarray}
This establishes the first condition for the validity of the Mountain Pass Theorem, for the functional $\mathbf{E}$.
To check the behaviour of $\mathbf{E}$ at infinity, let us consider $\phi=t e$ where $t\in\mathbb{R}$ and $e \in \mathbb{C}^{2m+1}$, $\mid\mid e\mid\mid_{\ell^2_{\epsilon}}=1$. We have that
\begin{eqnarray}
\mathbf{E}(te)=\frac{t^2}{2}-\frac{1}{2\sigma +2}t^{2\sigma +2}||e||_{\ell^{2\sigma +2}}^{2\sigma +2} \rightarrow -\infty ,
\nonumber
\end{eqnarray}
as $t\rightarrow +\infty$. This ensures the validity of the second condition of the Mountain Pass Theorem. For fixed $\phi\neq 0$ and choosing $t$ sufficiently large, we may set $u_1=t\phi$ to obtain the existence of a non-trivial solution for  (\ref{swef})-(\ref{swefc}).\ \ $\diamond$
\begin{remark} {\em In theorem \ref{swee} we have to state explicitly that we deal with the finite dimensional problem, since otherwise we may not prove the validity of the PS condition. In the infinite dimensional case, due to the lack of the Schur property for the space $\ell^{2}$ (in contrast with the space $\ell^1$ which posses this property-weak convergence coincides with strong convergence), we may not conclude the strong convergence of the sequence, from its weak convergence. The generalization of the mountain pass argument in the case of the infinite dimensional conservative system is a subject of future research, and  will be considered elsewhere. }\  $\bullet$
\end{remark} 
%%%%%%%%%%%%%%%%%%%%%%%%%%%%%%%%%%%%%%%%%%%%%%%%%%%%%%%
%%%%%%%%%%%%%%%%%%%%%%%%%%%%%%%%%%%%%%%%%%%%%%%%%%%%%%%
%%%%%%%Nonexistence%%%%%%%%%%%%%%%%%%%%%%%%%%%%%%%%%%%
\subsection{Non-existence of non trivial standing waves}
We now provide some non-existence results for nontrivial standing waves. The results that follow are valid for both problems (\ref{swef})-(\ref{swefc}) and (\ref{swe}).

We recall \cite[Theorem 18.E, pg. 68]{zei85}
(Theorem of Lax and Milgram), which will be used to establish existence of solutions for an auxiliary linear, non-homogeneous problem, related to (\ref{swe}). 
\begin{rigor1}
\label{LMth} 
Let $X$ be a Hilbert space and $A:X\rightarrow X$ be a linear continuous operator. Suppose that there exists  $c^*>0$ such that
\begin{eqnarray}
\label{strongmonot}
\mathrm{Re}(Au,u)_X\geq c^*||u||^2_X,\;\;\mbox{for all}\;\;u\in X.
\end{eqnarray}
Then for given $f\in X$, the operator equation $Au=f,\;\;u\in X$, has a unique solution
\end{rigor1}
\begin{rigor1}
There exist no nontrivial standing waves  of energy less than $E_c(\omega,\sigma):=\left(\frac{\omega^4}{4}\right)^{1/4\sigma}$.
\end{rigor1}
{\bf Proof:}\ \  Let $\omega\in\mathbb{R}$, $\omega\neq 0$, be fixed. For any $\epsilon >0$, the operator $A_{\omega}:\ell^2\rightarrow\ell^2$, defined by
\begin{eqnarray}
\label{strongop1}
(A_{\omega}\phi)_{n\in\mathbb{Z}}&=&-\frac{1}{\epsilon}(\phi_{n-1}-2\phi_n+\phi_{n+1})+\omega^2\phi_n,
\end{eqnarray}
is linear and continuous and satisfies assumption (\ref{strongmonot}) of Theorem \ref{LMth}: Using (\ref{lpno}) and (\ref{lp4}) we get that
\begin{eqnarray}
\label{check}
(A_{\omega}\phi,\phi)_{\ell^2}=\frac{1}{\epsilon}||B\phi||^2_{\ell^2}+\omega^2||\phi||^2\geq \omega^2||\phi||^2_{\ell^2}\;\;\mbox{for all}\;\;\phi\in\ell^2.
\end{eqnarray}
For given $z\in\ell^2$, we consider the linear operator equation
\begin{eqnarray}
\label{linear}
(A_{\omega}\phi)_{n\in\mathbb{Z}}=|z_n|^{2\sigma}z_n,\;\;\phi\in\ell^2.
\end{eqnarray}
The function $F_1:\mathbb{C}\rightarrow\mathbb{C}$, $F_1(z)=|z_n|^{2\sigma}z_n$ satisfies Lemma \ref{LipN}, with $\beta=2\sigma$. Therefore assumptions of Theorem \ref{LMth} are satisfied, and (\ref{linear}) has a unique solution $\phi\in\ell^2$. For some $R>0$, we consider the closed ball of $\ell^2$, $B_R:=\{z\in\ell^2\;:||z||_{\ell^2}\leq R\}$, and we define the map
$\mathcal{P}:\ell^2\rightarrow\ell^2$, by $\mathcal{P}(z):=\phi$ where $\phi$ is the unique solution of the operator equation (\ref{linear}). Clearly the map $\mathcal{P}$ is well defined. 

Let $z$, $\xi\in B_R$ such that $\phi=\mathcal{P}(z)$, $\psi=\mathcal{P}(\xi)$. The difference $\chi:=\phi-\psi$ satisfies the equation
\begin{eqnarray}
\label{claim2} 
(A_{\omega}\chi)_{n\in\mathbb{Z}}=F_1(z)-F_1(\xi).
\end{eqnarray}
Taking now the scalar product of (\ref{claim2}) with $\chi$ and using (\ref{ufp2}), we derive 
\begin{eqnarray}
\label{cmap1a}
\frac{1}{\epsilon}||B\chi||^2_{\ell^2}+\omega^2||\chi||^2_{\ell^2}&\leq& ||F_1(z)-F_1(\xi)||_{\ell^2}||\chi||_{\ell^2}\nonumber\\
&\leq&2R^{2\sigma}||z-\xi||_{\ell^2}||\chi||_{\ell^2}\nonumber\\
&\leq&\frac{\omega^2}{2}||\chi||_{\ell^2}^2+\frac{2}{\omega^2}R^{4\sigma}||z-\xi||_{\ell^2}^2.
\end{eqnarray}
From (\ref{cmap1a}), (\ref{lp4}), (\ref{lpno}), we obtain the inequality
\begin{eqnarray}
\label{claim3}
||\chi||_{\ell^2}^2=||\mathcal{P}(z)-\mathcal{P}(\xi)||_{\ell^2}^2\leq \frac{4}{\omega^4}R^{4\sigma}||z-\xi||^2_{\ell^2}.
\end{eqnarray}
Since $\mathcal{P}(0)=0$, from inequality  (\ref{claim3}), we derive that for $R< (\omega^4/4)^{1/4\sigma}:=E_c(\omega,\sigma)$,  the map  $\mathcal{P}:B_R\rightarrow B_R$ and is a contraction.   Therefore  $\mathcal{P}$, satisfies the assumptions of Banach Fixed Point Theorem and has a unique fixed point, the trivial one. Hence, for $R<E_c(\omega,\sigma)$ the only standing wave is the trivial. \ \ $\diamond$   
\begin{remark}
{\em The above theorem is interesting when viewed in the following way: The $\ell^2$ norm
corresponds to the energy of the excitations of the lattice. If the energy of the excitation
is less that $E_c(\omega,\sigma)$ the lattice may not support a standing wave of frequency $\omega$.  This relation may be seen as some kind of dispersion relation of frequency vs energy for the standing wave solutions. It contains information on the type of nonlinearity, through its dependence on the nonlinearity exponent $\sigma$. Let us note that a similar computation to (\ref{cmap1a}) shows that  $\mathcal{P}$ maps also the ball $B_R$  with $R\leq (\omega^4)^{1/4\sigma}:=R_1$, to itself and obviously $E_c<R_1$, but wthin this range, the map is not essentially a contraction.} \ \ $\bullet$ 
\end{remark} 
\begin{remark}
{\em Arguments of the same type as above, may be used for the proof of existence of unique steady states for the case of nonlinearities that do not vanish when $\phi=0$ (see the end of Remark \ref{ENERGY})}.\ \ $\bullet$
\end{remark} 
  
%%%%%%%%%%%%%%%%%%%%%%%%%%%%%%%%%%%%%%%%%%%%%%%%%%%%%%%%%%%%%%%%%%%%%%%%%%%%%%%%%%%%%%%%%%%
%%%%%%%%%%%%%%%%%%%%%%%%%%%%%%%%%%%%%%%%%%%%%%%%%%%%%%%%%%%%%%%%%%%%%%%%%%%%%%%%%%%%%%%%%%%
\section{Existence of global attractor in $\ell^2$ for the dissipative case}
In this section we study the asymptotic behavior of solutions of (\ref{lat1})-(\ref{lat2}) in the dissipative case, $\delta >0$ and $g_n\neq 0$, i.e. taking into account the effect of a weak dissipation and of an external excitation. This time, as a
model problem we shall employ the dissipative case of (\ref{cons1})-(\ref{cons2})
\begin{eqnarray}
\label{dissip1}
i\dot{u}_n+\frac{1}{\epsilon}(u_{n-1}-2u_n+u_{n+1})&+&i\delta u_n +|u_n|^{2\sigma}u_n=g_n,\;\;n\in\mathbb{Z},\;\;0\leq\sigma <\infty,\;\;\epsilon>0,\\
\label{dissip2}
u_n(0)&=&u_{n,0},\;\;n\in\mathbb{Z}.
\end{eqnarray}
We also consider
the weakly damped driven NLS partial differential equation
\begin{eqnarray}
\label{dpde}
i\partial_t u + u_{xx}+i\delta u&+&|u|^{2\sigma}u=g,\;\; x\in\mathbb{R},\;\;t>0,\\
u(x,0)&=&u_0(x).\nonumber
\end{eqnarray}
It is well known that the dissipation and forcing terms have
important effects in the long-time behavior of (\ref{dpde}). As an
example, we refer to the case $\sigma=2$,  for which has been
observed numerically the existence of chaotic attractors
\cite{blow,Nozaki}, in constrast with the conservative case which
is completely integrable  by
the inverse scattering theory, for sufficiently smooth initial data. We refer to
\cite{chid88,ogoub,nik02,laur95,XWang} for results on the
existence, finite dimensionality, and regularity of global
attractors for (\ref{dpde}).
\subsection{Existence of absorbing set}
This subsection is devoted to the existence of an absorbing set in
$\ell^2$. We have the following Lemma.
\begin{smallTheorem}
\label{exabsball}
Let $(u_{0,n})_{n\in\mathbb{Z}}=u_0, (g_n)_{n\in\mathbb{Z}}=g\in\ell^2$. For any $\epsilon >0$ and $0\leq\sigma <\infty$, the dynamical system defined by
(\ref{dissip1})-(\ref{dissip2})
\begin{eqnarray}
\label{dynamical} S(t):u_0\in {\ell}^2\rightarrow u(t)\in
{\ell}^2,
\end{eqnarray}
possesses a bounded absorbing set $\mathcal{B}_0$ in ${\ell}^2$:
For every bounded set $\mathcal{B}$ of ${\ell}^2$, there exists
$t_0(\mathcal{B},\mathcal{B}_0)$ such that for all  $t\geq
t_0(\mathcal{B},\mathcal{B}_0)$, it holds
$S(t)\mathcal{B}\subset\mathcal{B}_0$.
\end{smallTheorem}
{\bf Proof:}\ \. We take the scalar product of (\ref{dissip1}) with $iu(t)$ and we get
\begin{eqnarray}
\label{eneg0}
\frac{1}{2}\frac{d}{dt}||u||^2_{\ell^2}+\delta||u||^2_{\ell^2}=\mathrm{Im}\sum_{n\in\mathbb{Z}}\overline{u}_ng_n\leq\frac{\delta}{2}||u||^2_{\ell^2}+\frac{1}{2\delta}||g||_{\ell^2}^2.
\end{eqnarray}
and therefore
\begin{eqnarray}
\label{absset1}
\frac{d}{dt}||u||^2_{\ell^2}+\delta||u||^2_{\ell^2}\leq \frac{1}{\delta}||g||_{\ell^2}^2.
\end{eqnarray}
From (\ref{absset1}) we derive that
$u\in\mathrm{L^{\infty}}([0,\infty),\ell^2)$: Gronwall's Lemma
implies that
\begin{eqnarray}
\label{gronw}
||u(t)||_{\ell^2}^2\leq ||u_0||^2_{\ell^2}\exp(-\delta t)+\frac{1}{\delta^2}||g||_{\ell^2}^2
\{1-\exp(-\delta t)\}.
\end{eqnarray}
Letting $t\rightarrow\infty$ we infer that
\begin{eqnarray*}
\label{abbset2}
\limsup_{t\rightarrow\infty}||u(t)||^2_{\ell^2}\leq\frac{1}{\delta^2}||g||_{\ell^2}^2.
\end{eqnarray*}
Setting $\rho^2=||g||^2_{\ell^2}/\delta^2$, it follows that for any number $\rho_1>\rho$ the ball $\mathcal{B}_0$ of $\ell^2$ centered at $0$ of radius $\rho_1$ is an absorbing set
for the semigroup $S(t)$: If $\mathcal{B}$ is a bounded set of $\ell^2$ included in a ball of $\ell^2$ centered at $0$ of radius $R$, then for 
$t\geq t_0(\mathcal{B},\mathcal{B}_0)$ where
\begin{eqnarray}
\label{time}
t_0=\frac{1}{\delta}\log\frac{R^2}{\rho_1^2-\rho^2},
\end{eqnarray}
it holds $||u(t)||^2_{\ell^2}\leq\rho_1^2$, i.e. $S(t)\mathcal{B}\subset\mathcal{B}_0$.
Note that in the absence of external excitation, the dynamical system exhibits trivial dynamics, in the sense that  $\limsup_{t\rightarrow\infty}||u(t)||_{\ell^2}^2=0$, as (\ref{gronw}) clearly shows.\ \ $\diamond$
\begin{remark}
\label{Absl1}
{\em The existence of the absorbing set, can be shown again by working in a different way, in the equivalent norm of $\ell^2_1$:
By taking the scalar product of (\ref{dissip1}) with $-\dot{u}-\delta u$ we may derive 
the energy equation
\begin{eqnarray}
\label{enegs1}
\frac{1}{2}\frac{d}{dt}\mathrm{J}(u(t))&+&\delta\mathrm{J}(u(t))=\Lambda (u(t)),\;\;
\mathrm{J}(u):=||u||_{\ell^2_1}^2-\epsilon\left\{\frac{1}{\sigma +1}||u||_{\ell^{2\sigma +2}}^{2\sigma +2}-2(g,u)\right\},\\
\Lambda (u)&:=&\epsilon\delta\left\{\frac{\sigma}{\sigma +1}
||u||_{\ell^{2\sigma +2}}^{2\sigma +2}+(g,u)_{\ell^2}\right\}
+\mathrm{Im}\sum_{n\in\mathbb{Z}}\overline{u}_ng_n.
\end{eqnarray}
From (\ref{lp1}) and (\ref{time}) we obtain the following estimate for $\Lambda (u)$
\begin{eqnarray}
\label{enegs2}
\Lambda (u)\leq \epsilon\delta\left\{\frac{\sigma}{\sigma +1}\rho_1^{2\sigma +2}+ ||g||_{\ell^2}\rho_1\right\}+||g||_{\ell^2}\rho_1=\lambda_1.
\end{eqnarray}
Gronwall's inequality applied in (\ref{enegs1}) implies that
\begin{eqnarray}
\mathrm{J}(u(t))\leq\mathrm{J}(u_0)\exp (-2\delta t)+\frac{\lambda_1}{\delta}
\left\{1-\exp(-2\delta t)\right\}
\end{eqnarray}
Letting $t\rightarrow\infty$ we get $\limsup_{t\rightarrow\infty}\mathrm{J}(u(t))\leq \lambda_1/\delta$. Hence, for $t\geq t_0$
\begin{eqnarray}
\label{enegs3}
||u(t)||^2_{\ell^2_1}\leq \epsilon\rho_1^{2\sigma +2}
+3\epsilon ||g||_{\ell^2}\rho_1+\frac{1}{\delta}||g||_{\ell^2}\rho_1:=\rho_2^2
\end{eqnarray}
Relation (\ref{enegs3}), gives an estimate of the radius of the absorbing ball, involving the degree of the nonlinearity and lattice spacing. Note for example, that for sufficiently small $$0<\epsilon<\frac{\rho_1-||g||_{\ell^2}\delta^{-1}}{\rho_1^{2\sigma+1}+3||g||_{\ell^2}},$$ we get that $\rho_2<\rho_1$. }$\bullet$
\end{remark}
Throughout the rest of the paper, for simplicity reasons we set $\epsilon=1$.
%%%%%%%%%%%%%%%%%%%%%%%%%%%%%%%%%%%%%%%%%%%%%%%%%%%%%%%%%%%%%%%%%%%%%%
%%%%%%%%%%%%%%%%%%%%%%%%%%%%%%%%%%%%%%%%%%%%%%%%%%%%%%%%%%%%%%%%%%%%%%
%%%%%%%%%%%%%%%%%%%%%%%%%%%%%%%%%%%%%%%%%%%%%%%%%%%%%%%%%%%%%%%%%%%%%%
\subsection{Asymptotic compactness of the semigroup and existence of global attractor}
We now prove the asymptotic compactness of the semigroup generated by the DNLS equation ({\ref{dissip1})-(\ref{dissip2}).
To establish this property, we follow the approach introduced in \cite{Bates} and applied also in \cite{SZ1} for a second order lattice dynamical system. This approach is based on the derivation of estimates on the tail ends of solutions of ({\ref{dissip1})-(\ref{dissip2}). 
\begin{smallTheorem}
\label{DNLStail}
Let $(u_{0,n})_{n\in\mathbb{Z}}=u_0\in\mathcal{B}$ where $\mathcal{B}$ is a bounded set of $\ell^2$ and $(g_n)_{n\in\mathbb{Z}}=g\in\ell^2$. For any $0\leq\sigma <\infty$, and $\eta >0$, there exist $T(\eta )$ and $K(\eta )$ such that the solution
$(u_{n})_{n\in\mathbb{Z}}=u$ of (\ref{dissip1})-(\ref{dissip2}) satisfies  for all $t\geq T(\eta )$, the estimate
\begin{eqnarray}
\label{precom}
\sum_{\mid n\mid >2M}\mid u_n(t)\mid^2 \le \frac{2\eta}{\delta},\;\;\mbox{for any}\;\;
M>K(\eta ).
\end{eqnarray}
\end{smallTheorem}
{\bf Proof:}\ Choose a smooth function $\theta \in C^{1}({\mathbb R}^{+},{\mathbb R})$ satisfying the following properties
\begin{eqnarray}
\left\{
\begin{array}{ccc}
\theta(s)=0,\;\; & 0\le s\le 1  \nonumber \\
0\le \theta(s) \le 1,\;\; & 1\le s\le 2 \nonumber \\
\theta(s)=1,\;\; & s\ge 2. \nonumber
\end{array}
\right.
\end{eqnarray}
and
\begin{eqnarray}
\label{mvth}
\mid \theta^{'}(s)\mid \le C_0, \hspace{2mm} s\in {\mathbb R}^{+},
\end{eqnarray}
for some $C_0 \in {\mathbb R}$. 
We shall use the shorthand notation $\theta_n=\theta\left( \frac{\mid n\mid}{M}\right)$ and $\mathrm{Re}(u_n)=u_{1,n}$, $\mathrm{Im}(u_n)=u_{2,n}$
We now multiply the DNLS equation with the function $\theta_n\bar{u}_n$, $n\in {\mathbb Z}$, and we sum over all sites and keep the imaginary part.
The effect of the  function $\theta$ on $u_n$ is to cut-off all sites within a ball of radius comparable to $M$ and to take into account only the remote sites.
This results to
\begin{eqnarray}
\frac{1}{2}\frac{d}{dt}\sum_{n\in\mathbb{Z}}\theta_n\mid u_n\mid^2 +
\delta\sum_{n\in\mathbb{Z}}\theta_n\mid u_n\mid^2
-\sum_{n\in\mathbb{Z}} \{ (Bu_{2})_n(B\theta
u_{1})_n-(Bu_{1})_n(B\theta u_{2})_n\}
= \mathrm{Im}\sum_{n\in\mathbb{Z}}g_n \theta_n\bar{u}_n . \label{eqn:dNLStail}
\end{eqnarray}
The nonlinear term contributes a purely real term so it drops out of the above equality.
We must now estimate the remaining terms.
Using the Cauchy-Schwartz inequality we find for the term containing the forcing that
\begin{eqnarray}
\mathrm{Im}\sum_{n\in\mathbb{Z}}g_n \theta_n\bar{u}_n \nonumber \le
\frac{\delta}{2}\sum_{n\in\mathbb{Z}} \theta_n\mid u_n\mid^2 +\frac{1}{2\delta} \sum_{\mid
n\mid \ge M}\mid g_n\mid^2 . \nonumber
\end{eqnarray}
We now estimate the last term on the lhs of  equation (\ref{eqn:dNLStail}). After some algebra we obtain
\begin{eqnarray}
I_1=\sum_{n\in\mathbb{Z}} \{ (Bu_{2})_n(B\theta
u_{1})_n-(Bu_{1})_n(B\theta u_{2})_n\}=
\sum_{n\in\mathbb{Z}}\{(\theta_{n+1}-\theta_{n})(u_{1,n}u_{2,n+1}-u_{2,n}u_{1,n+1})\}.
\nonumber
\end{eqnarray}
We may further estimate $I_1$ as follows,
\begin{eqnarray}
|I_1|&\leq& \sup_{n\in\mathbb{Z}}\mid \theta_{n+1}-\theta_{n}\mid\; \sum_{n \in {\mathbb Z}} \left|(u_{1,n}u_{2,n+1}-u_{2,n}u_{1,n+1})\right|
\nonumber \\
&\leq& \frac{C_1}{M}  \mid\mid u\mid\mid_{\ell^{2}}^2 \leq \frac{C_1}{M}\rho_1^2,
\nonumber
\end{eqnarray}
where $C_1$ depends only on $C_0$. For the last inequality we applied the mean value theorem and (\ref{mvth}) 
and the restriction of the solution to the absorbing ball in $\ell^{2}$ ($\rho_1$ is the radius of the
absorbing ball of Lemma \ref{exabsball}).

Therefore we obtain the following differential inequality
\begin{eqnarray}
\frac{1}{2}\frac{d}{dt}\sum_{n\in\mathbb{Z}} \theta_n \mid
u_n\mid^2  +\frac{\delta}{2}\sum_{n\in\mathbb{Z}}\theta_n \mid
u_n\mid^2 \le \frac{C_1}{M}\rho_1^2 +\frac{1}{2\delta}\sum_{\mid
n \mid
>M}\mid g_n\mid^2 .\nonumber
\end{eqnarray}
Using the Gronwall inequality we obtain the following  estimate
\begin{eqnarray}
\sum_{n\in\mathbb{Z}}\theta_n \mid u_n \mid^2 \le
e^{-\delta(t-t_0)}\sum_{n\in\mathbb{Z}}\theta_n \mid
u_n(t_0)\mid^2 +\frac{1}{\delta}\left( \frac{2C_1}{M}\rho_1^2
+\frac{1}{\delta}\sum_{\mid n \mid >M}\mid g_n\mid^2  \right),
\nonumber
\end{eqnarray}
for $t>t_0$ where $t_0$ is the time of entry of initial data bounded in $\ell^2$, to the absorbing ball of radius $\rho_1$ in $\ell^{2}$.
Since $g\in \ell^{2}$, then  for all $\eta >0$, there exists $K(\eta)$ such that
\begin{eqnarray}
\frac{2C_1}{M}\rho_1^2 +\frac{1}{\delta}\sum_{\mid n \mid >M}\mid g_n\mid^2  \le \eta, \hspace{3mm} \forall M > K(\eta).
\nonumber
\end{eqnarray}
Therefore, for all $\eta$ and for $t>t_0$ and $M>K(\eta)$ we obtain
that
\begin{eqnarray}
\sum_{n\in\mathbb{Z}}\theta_n \mid u_n \mid^2 \le
e^{-\delta(t-t_0)}\rho_1^2 +\frac{1}{\delta}\eta \nonumber
\end{eqnarray}
Choosing $t$ large enough we may then obtain
\begin{eqnarray*}
\sum_{\mid n\mid >2M}\mid u_n\mid^2 \le
\sum_{n\in\mathbb{Z}}\theta_n\mid u_n\mid^2 \le
\frac{2\eta}{\delta} .\nonumber
\end{eqnarray*}
This estimate holds as long as $t\ge T(\eta)$ where
$$T(\eta)=t_0+\frac{1}{\delta}\ln\left(\frac{\delta\rho_1^2}{\eta}\right)$$
and $M> K(\eta)$. This concludes the proof of the Lemma.
$\diamond$
\\ \\
{\bf Remark}  The above
Lemma has an additional implication, since  it allows us to
obtain an estimate for the time scales  required for the
localization of the energy to the part of the lattice with $\mid
n\mid\leq M$. \ $\bullet$
\begin{rigor2}
\label{asymcomp}
Let $0\leq\sigma <\infty$. The semigroup $S(t)$ is asymptotically compact in $\ell^{2}$, that is, if the sequence $\phi_n$ is bounded in $\ell^{2}$ and $t_n \rightarrow \infty$, then $S(t_n)\phi_n$ is precompact in $\ell^2$.
\end{rigor2}
{\bf Proof:} Let us consider a bounded sequence $\phi_n \in
\ell^2$, such that $\mid\mid \phi_n\mid\mid_{\ell^2}\le r$, $r>0$,
$n=1,2,...$ By Lemma \ref{exabsball}, there exists $T(r)>0$ such
that $S(t)\phi_n \subset {\cal B}_0$, $\forall t\ge T(r)$. As
$t_n\rightarrow \infty$ there exists an integer $N_1(r)$ such that
$t_n \ge T(r)$ for $n \ge N_1(r)$ and
\begin{eqnarray}
\label{prese}
S(t_n)\phi_n \subset {\cal B}_0, \hspace{3mm} \forall n \ge N_1(r).
\end{eqnarray}
Therefore $S(t_n)\phi_n$ is weakly relatively compact. Thus, there exists $\phi_0 \in \mathcal{B}_0$ and a subsequence of $S(t_n)\phi_n$ (not relabelled) such that $S(t_n)\phi_n \rightarrow \phi_0$ weakly in $\mathcal{B}_0$. To ensure  precompactness we need to show that the above convergence is strong, that is for all $\eta >0$, there exists $N(\eta) \in {\mathbb N}$ such that
$$ \mid\mid S(t_n)\phi_n - \phi_0 \mid\mid_{\ell^{2}} \le \eta, \hspace{3mm} \forall n \ge N(\eta).$$
From (\ref{prese}) and the tail estimates provided by Lemma
\ref{DNLStail}, we may conclude the existence of some $K_1(\eta)$
and $T_1(\eta)$ such that
\begin{eqnarray}
\sum_{\mid i\mid \ge K_1(\eta)} \mid (S(t)S(T_r)\phi_n)_i \mid^{2} \le \frac{\eta^2}{8}, \hspace{3mm} t \ge T_1(\eta),
\nonumber
\end{eqnarray}
where by $(S(t)\phi_n )_i$ we denote the $i$-th coordinate of the infinite sequence $S(t)\phi_n \in \ell^{2}$.
Since $t_n \rightarrow \infty$ we may find $N_2(r,\eta)\in {\mathbb N}$ such that $t_n \ge T_r +T(\eta)$ if $n \ge N_2(r,\eta)$. Hence,
\begin{eqnarray}
\sum_{\mid i\mid \ge K_1(\eta)} \mid (S(t_n)\phi_n)_i \mid^{2} = \sum_{\mid i\mid \ge K_1(\eta)}\mid (S(t_n-T_r)S(T_r)\phi_n)_i \mid^{2} \le \frac{\eta^2}{8}
\nonumber
\end{eqnarray}
where we used the semigroup property for $S(t)$.

Since $\phi_0\in \ell^{2}$ we have the following tail estimate
$$ \sum_{\mid i\mid \ge K_2(\eta)} \mid (\phi_0)_i \mid^2 \le \frac{\eta^2}{8}.$$
Let us now choose $K(\eta)=max(K_1(\eta),K_2(\eta))$. Since $S(t_n)\phi_n \rightharpoonup  \phi_0$ in $\ell^2$, we may see that $(S(t_n)\phi_n)_{i} \rightarrow (\phi_0)_{i}$ strongly in ${\mathbb C}^{2K(\eta)+1}$. To prove the strong convergence we will break $\mid\mid \cdot \mid\mid_{\ell^{2}}$ into two distinct parts, the finite dimensional part (corresponding to ${\mathbb C}^{2K(\eta)+1}$) and the infinite dimensional part corresponding to the tails. Thus,
\begin{eqnarray}
\mid\mid S(t_n)\phi_n -\phi_0 \mid\mid_{\ell^{2}}^2 = \sum_{\mid i\mid \le K(\eta)}\mid (S(t_n)\phi_n -\phi_0)_i\mid^2 +
\sum_{\mid i\mid > K(\eta)}\mid (S(t_n)\phi_n -\phi_0)_i\mid^2
\nonumber \\
\le \frac{\eta^2}{2}+ 2 \sum_{\mid i\mid > K(\eta)}(\mid (S(t_n)\phi_n)_i \mid^2 + \mid(\phi_0)_i\mid^2) \le \eta^{2}
\nonumber
\end{eqnarray}
where the first estimate comes from the strong convergence in the finite dimensional space ${\mathbb C}^{2K(\eta)+1}$. \ $\diamond$

Using  Proposition \ref{asymcomp} and Theorem 1.1.1 of
\cite{RTem88}, we may arrive at the main result of this section
which can be stated as follows:
\begin{rigor1}
\label{dynamics1} Let $0\leq\sigma <\infty$. 
The semigroup $S(t)$ associated to (\ref{dissip1})-(\ref{dissip2})
possesses a global attractor
$\mathcal{A}=\omega(\mathcal{B}_0)\subset\mathcal{B}_0\subset\ell^2$
which is compact, connected and maximal among the functional
invariant sets in $\ell^2$.
\end{rigor1}
%%%%%%%%%%%%%%%%%
%%%%%%%%%%%%%%%%%
%%%%%%%%%%%%%%%%%
\section{Finite dimensional   approximation  of the global attractor}
This section is devoted to a result of stability of the global attractor of the dynamical system generated by (\ref{dissip1})-(\ref{dissip2}) under its approximation by a global attractor of an appropriate finite dimensional dynamical system.

Let us note first that problem (\ref{dissip1})-(\ref{dissip2}) could be viewed as a boundary value problem satisfying the  boundary condition that the solutions are bounded at infinity. Since $u(t)\in\mathrm{C}^1(\mathbb{R}^+,\ell^2)$ clearly $$\lim_{n\rightarrow\infty}u_n(t)=0,\;\;t\geq0.$$
We seek for an approximation of (\ref{dissip1})-(\ref{dissip2})
by a finite dimesional system of ordinary differential equations.  

We consider the boundary value problem in $\mathbb{C}^{2m+1}$,
\begin{eqnarray}
\label{bvp1}
i\dot{\psi}_n+(\psi_{n-1}-2\psi_n+\psi_{n+1})&+&i\delta \psi_n +|\psi_n|^{2\sigma}\psi_n=g_n,\;\;|n|\leq m,\;\;0\leq\sigma <\infty,\\
\label{bvp2}
\psi_{-(m+1)}(\cdot)&=&\psi_{(m+1)}(\cdot)=0,\\
\label{bvp3}
\psi_n(0)&=&\psi_{n,0},\;\;|n|\leq m.
\end{eqnarray}
In a similar manner to the infinite dimensional system we may show that for the well posedness and asymptotic behavior of solutions of (\ref{bvp1})  the following result holds.
\begin{rigor2}
\label{fattr}
Let $\psi_0:=(\psi_{n,0})_{|n|\leq m}\in\mathbb{C}^{2m+1}$. For $0\leq\sigma <\infty$, there exists a unique solution of
(\ref{bvp1})-(\ref{bvp3})
such that $\psi\in\mathrm{C}^1([0,\infty),\mathbb{C}^{2m+1})$. The dynamical system defined by  (\ref{bvp1})-(\ref{bvp3})
\begin{eqnarray}
\label{fdynamical} S_m(t):\psi_0\in \mathbb{C}^{2m+1}\rightarrow
\psi(t)\in \mathbb{C}^{2m+1},
\end{eqnarray}
possesses a bounded absorbing set $\mathcal{O}_0$ in
$\mathbb{C}^{2m+1}$ and a global attractor
$\mathcal{A}_{m}\subset\mathcal{O}_0\subset\mathbf{C}^{2m+1}$ :
For every bounded set $\mathcal{O}$ of ${\ell}^2$, there exists
$t_1(\mathcal{O},\mathcal{O}_0)$ such that for all  $t\geq
t_1(\mathcal{O},\mathcal{O}_0)$, it holds that
$S_m(t)\mathcal{O}\subset\mathcal{O}_0$ and for every $t\geq 0$
$S_m(t)\mathcal{A}_m=\mathcal{A}_m$.
\end{rigor2}
Note that the estimates used for the proof of Proposition
\ref{fattr} (similar to those used in the proof Lemma
\ref{exabsball} {\em are independent of $m$}. 

Following \cite{SZ1,SZ2} (see also
\cite{bab90,2Nikoi-a} for a similar idea applied to pdes
considered in all of $\mathbb{R}^N$), we observe that the
$\mathbb{C}^{2m+1}$-solution of (\ref{bvp1})-(\ref{bvp3}) can be
extended naturally in the infinite dimensional space $\ell^2$ , as
\begin{equation}
\label{ZN}
(\phi_m(t))_{m\in\mathbb{Z}}:=\left\{
\begin{array}{cc}
\psi(t):=(\psi_n(t))_{|n|\leq m},\;\; & |n|\leq m,\\
0,\;\; & |n|> m. \nonumber
\end{array}
\right.
\end{equation}

The extension (\ref{ZN}) will be used in order to verify that the global attractor
$\mathcal{A}$ of the semigroup $S(t)$ associated with
(\ref{dissip1})-(\ref{dissip2}), can be approximated by the global
attractor $\mathcal{A}_m$ of $S_m(t)$ associated to
(\ref{bvp1})-(\ref{bvp3}) as $m\rightarrow\infty$. We recall that
the semidistance of two nonempty compact subsets of a metric space
$X$ endowed with the metric $d_X(\cdot ,\cdot)$ is defined as
\begin{eqnarray*}
d(\mathcal{B}_1,\mathcal{B}_2)=\sup_{x\in\mathcal{B}_1}\inf_{y\in\mathcal{B}_2}d_X(x,y).
\end{eqnarray*}
\begin{rigor1}
\label{dynamics2}
For any  $0\leq\sigma <\infty$, the global attractor $\mathcal{A}_m$ converges to $\mathcal{A}$ in the sense of the semidistance related to $\ell^2$: $\lim_{m\rightarrow\infty}d(\mathcal{A}_m,\mathcal{A})=0$.
\end{rigor1}
{\bf Proof:}\ \ We denote by $\mathcal{U}$ an open-neighborhood of
the absorbing ball $\mathcal{B}_0$ of $S(t)$. Obviously
$\mathcal{A}$ attracts $\mathcal{U}$. For arbitrary
$m\in\mathbb{Z}$, we consider the  semigroup $S_m(t)$ defined by
Proposition \ref{fattr} and its global attractor $\mathcal{A}_m$.
Exactly as in Lemma \ref{exabsball}, it can be shown that
$\mathcal{B}_0\cap\mathbb{C}^{2m+1}$ is also an absorbing set
for $S_m(t)$. Therefore
$$\mathcal{A}_m\subset\mathcal{B}_0\cap\mathbb{C}^{2m+1}\subset\mathcal{U}\cap\mathbb{C}^{2m+1},$$
and $\mathcal{A}_m$ attracts $\mathcal{U}\cap\mathbb{C}^{2m+1}$.
In the light of Proposition \ref{fattr} and \cite[Theorem
I1.2, pg. 28]{RTem88}, it remains to verify that for every compact interval
$\mathrm{I}$ of $\mathbb{R}^+$,
\begin{eqnarray}
\label{cond2} \delta_m(\mathrm{I}):=\sup_{\psi_0\in\mathcal{U}\cap\mathbb{C}^{2m+1}}
\sup_{t\in\mathrm{I}}d(S_m(t)\psi_0,S(t)\psi_0)\rightarrow
0,\;\;\mbox{as}\;\;m\rightarrow\infty .
\end{eqnarray}
We consider the corresponding solution $\psi(t)=S_m(t)\psi_0$,
$\psi(0)=\psi_0$, in $\mathbb{C}^{2m+1}$ through
(\ref{bvp1})-(\ref{bvp3}). Then by Proposition \ref{fattr},  it
follows that $\psi(t)\in\mathcal{A}_m$ for any $t\in\mathbb{R}^+$.
Therefore, if $\rho>0$ is the $m$-independent radius of the
absorbing ball $\mathcal{O}_m$ in  $\mathbb{C}^{2m+1}$, then
for every $t\in\mathbb{R}^+$, $||\psi(t)||^2_2\leq\rho^2$.  
Using (\ref{ZN}), we clearly observe that $\phi_m(t)$
also satisfies the estimates
\begin{eqnarray}
\label{aprrox1}
||\phi_m(t)||_{\ell^2}^2\leq\rho^2,\;\;||\dot{\phi}_m(t)||_{\ell^2}^2\leq
C(\rho, ||g||),
\end{eqnarray}
the latter derived by (\ref{bvp1}).
According to \cite[Theorem 10.1 pg. 331-332]{RTem88} or \cite[Lemma 4, pg. 60]{SZ2}, for the  justification of (\ref{cond2}) it suffices to show that $\phi_m(t)$ converges to  a solution $\phi(t)$ of (\ref{dissip1})-(\ref{dissip2}) in an arbitrary compact interval of $\mathbb{R}^+$, and $\phi_0=\phi(0)$ in a bounded set of $\ell^2$. Let $\mathrm{I}$ be an arbitrary compact interval of $\mathbb{R}^+$. 
From estimates (\ref{aprrox1}), we may extract a subsequence $\phi_\mu$ of $\phi_m$, such that
\begin{eqnarray}
\label{dweakc} \phi_{\mu}(t)\rightharpoonup
\phi(t),\;\;\mbox{in}\;\;\ell^2,\;\;\mbox{as}\;\;
\mu\rightarrow\infty,\;\;\mbox{for every}\;\;t\in D,
\end{eqnarray}
where $D$ denotes a countable dense subset of $\mathrm{I}$.

For any $t\in \mathrm{I}$ we consider the  sequence
$$\gamma_m(t):=(\phi_m(t),z)_{\ell^2},\;\;z\in\ell^2,$$
which is diferrentiable by (\ref{aprrox1}) as a function of $t$ and $\gamma_m
'(t)=(\dot{\phi}_m(t),z)_{\ell^2}$. By the mean-value Theorem,
there exists $\xi\in \mathrm{I}$ such that, for fixed $t,s\in
\mathrm{I}$
\begin{eqnarray*}
|(\phi_m(t)-\phi_m(s),z)_{\ell^2}|=|\gamma_m(t)-\gamma_m(s)|&=&|(\dot{\psi}_m(t),z)_{\ell^2}|\,|t-s|\nonumber\\
&\leq&\sup_{\xi\in
I}||\dot{\phi}_m(\xi)||_{\ell^2}||z||_{\ell^2}|t-s|\leq C|t-s|.
\end{eqnarray*}
It follows then, that there exists an $m$-independent constant $C_1$ such that
\begin{eqnarray}
\label{Ascoli} ||\phi_m(t)-\phi_m(s)||_{\ell^2}\leq C_1|t-s|,
\end{eqnarray}
i.e the sequence $\gamma_m$ is equicontinuous.  Hence by Ascoli's Theorem, it follows that the convergence (\ref{dweakc}), holds uniformly on $\mathrm{I}$ as $\mu\rightarrow\infty$. Summarizing, we obtain that for the 
subsequence $\phi_\mu$ hold the convergence relations
\begin{eqnarray}
\label{passlim}
&&\phi_{\mu}\rightarrow\phi \;\;\mbox{in}\;\;\mathrm{C}(\mathrm{I},\ell^2),\nonumber\\
&&\phi_{\mu}\stackrel{*}{\rightharpoonup} \phi\;\;\mbox{in}\;\;\mathrm{L}^{\infty}(\mathrm{I},\ell^2),\\
&&\dot{\phi}_{\mu}\stackrel{*}{\rightharpoonup}
\dot{\phi}\;\;\mbox{in}\;\;\mathrm{L}^{\infty}(\mathrm{I},\ell^2).\nonumber
\end{eqnarray}
For every $z\in\ell^2$ and $\omega (t)\in C^{\infty}_0(\mathrm{I})$ we consider the formula (see \cite[p.g. 59]{cazS})
\begin{eqnarray}
\label{weakform}
\int_{\mathrm{I}}(i\dot{\phi}_{\mu}(t),z)_{\ell^2}\omega (t)dt -
\int_{\mathrm{I}}(B\phi_{\mu}(t),Bz)_{\ell^2}\omega (t)dt
&+&\delta\int_{\mathrm{I}}(i\phi_{\mu}(t),z)_{\ell^2}\omega (t)dt\nonumber\\
&+&\int_{\mathrm{I}}(f(\phi_{\mu}(t)),z)_{\ell^2}\omega (t)dt
=\int_{\mathrm{I}}(g,z)_{\ell^2}\omega (t)dt,
\end{eqnarray}
where $f(\phi_{\mu}):=|\phi_{\mu}|^{2\sigma}\phi_{\mu}$. Needless to say that any solution of (\ref{weakform}), is in our case,  a solution of (\ref{dissip1})-(\ref{dissip2}) and vice versa. Using
(\ref{passlim}) we may pass to the limit in (\ref{weakform}) .
Note for example, that since $f:\ell^2\rightarrow\ell^2$ is
Lipschitz continuous on bounded sets of $\ell^2$, there exists
from (\ref{aprrox1}), a constant $c(\rho)$ such that
$||f(\phi_{\mu})-f(\phi)||_{\ell^2}\leq
c(\rho)||\phi_{\mu}-\phi||_{\ell^2}$. Then from (\ref{passlim}) we
infer
\begin{eqnarray*}
\label{pasnon}
\left|\int_{\mathrm{I}}(f(\phi_{\mu}(t))-f(\phi(t)),z)\omega(t)dt\right|
&\leq&\int_{\mathrm{I}}||f(\phi_{\mu}(t))-f(\phi(t))||_{\ell^2}||z||_{\ell^2}\omega (t)dt\\
&\leq&c\int_{\mathrm{I}}||\phi_{\mu}(t)-\phi(t)||_{\ell^2}||z||_{\ell^2}\omega (t)\\
&\leq&c\sup_{t\in\mathrm{I}}||\phi_{\mu}(t)-\phi(t)||_{\ell^2}||z||_{\ell^2}\int_{\mathrm{I}}|\omega(t)|dt\rightarrow
0,\;\;\mbox{as}\;\;\mu\rightarrow\infty.
\end{eqnarray*}
Now we conclude by the same arguments as in \cite[Lemma 4, pg. 61]{SZ2}: Since $I$ is arbitrary, (\ref{weakform}) is satisfied for all $t\in\mathbb{R^+}$, i.e. $\phi(t)$ solves (\ref{dissip1})-(\ref{dissip2}). Moreover by (\ref{passlim}) we get that 
$\phi(t)$ is bounded in $\ell^2$ for all $t\in\mathbb{R^+}$. Therefore $\phi(t)\in\mathcal{A}$, which implies that $\phi_{\mu}(0)\rightarrow\phi(0)$ and $\phi_0$, is at least, in a bounded set of $\ell^2$.   Since the convergence holds for any other subsequence having the above formulated properties, by a contradiction argument using uniqueness, we may deduce that the convergence holds for the original sequence $\phi_m$. Condition
(\ref{cond2}) is proved.\ \ $\diamond$.
%%%%%%%%%%%%%%%%%%%%%
%%%%%%%%%%%%%%%%%%%%%
\section{The case of weighted spaces: Existence of spatially exponentially localized solutions}
We now turn to the properties of the solutions of the DNLS
(\ref{dissip1})-(\ref{dissip2}), in weighted spaces. We consider  
weight function $w_n$ which is an increasing function of $|n|$, satisfying for all $n\in\mathbb{Z}$, the following  condition:
\begin{eqnarray}
\mathrm{(W)} \left \{ \begin{array}{ccc}
&1\leq w_n\nonumber \\
&\mid w_{n+1}-w_{n} \mid \le d_1 w_n
\nonumber \\
&\underline{d}_2 w_n \le  w_{n+1}.
\nonumber
\end{array}
\right .
\end{eqnarray}
the space $\ell^2_w$, $$\ell_w^2=\{u_n \in {\mathbb C}:\;
\mid\mid u \mid\mid_{\ell_w^2}^2:=\sum_{n\in\mathbb{Z}} w_n \mid
u_n \mid^2 < \infty\}.$$  It can easily be seen that the space
$\ell_w^2$ is a Hilbert space with the norm $\mid\mid \;
\mid\mid_{\ell_w^2}$.  Such spaces are the discrete analogue of
weighted $L^{2}$ spaces. A choice for such a function may be the
exponential function $w_n=exp(\lambda \mid n\mid)$ for
$\lambda >0$. Existence of solutions in such spaces will provide
us with the existence of (exponentially) localized solutions for
the DNLS equation (\ref{dissip1})-(\ref{dissip2}). The use of such
spaces is important in the study of existence of soliton solutions
or breathers. An instance where such spaces have been used is in
\cite{Mackayexp} where  the existence of exponentially localized
solutions has been studied in conservative lattices using a
continuation argument, related to the anti-integrable limit.

Since the operator $A$ is not symmetric in the space $\ell_w^2$ we
cannot apply Theorem \ref{locex} for the local existence in such
spaces. We thus have to resort to general existence Theorems in
Banach spaces. To this end we need to show that the operators
involved are Lipschitz. This time as a model case we consider
(\ref{lat1})-(\ref{lat2}) with a nonlinearity satisfying
$\mathrm{(N_2)}$.

We have the following,
\begin{smallTheorem}\label{thm:LOCALLIPS}
Let condition $\mathrm{(N_2)}$ be fulfilled. The operator $T:\ell_{w}^{2}\rightarrow \ell_{w}^{2}$, defined by $T(z)=f(\mid z\mid^2)z$ satisfies the following properties:
\begin{itemize}
\item[(i)] $T$ is bounded on bounded sets of $\ell_{w}^{2}$ and
\item[(ii)] $T$ is locally Lipschitz continuous.
\end{itemize}
\end{smallTheorem}
{\bf Proof:} (i) Let $u\in B_R$, a closed ball of $\ell_{w}^2$ of
radius $R$. It follows from  condition $\mathrm{(W)}$ that $\mid z_n \mid^2
\le \mid\mid z \mid\mid_{\ell_{w}^{2}}^2$ for all $n \in {\mathbb
Z}$.  Since $f$ is
continuous, we may argue as for the proof of Lemma \ref{LipN}, in
order to get the inequality
\begin{eqnarray}
\mid\mid T(z) \mid\mid_{\ell_{w}^{2}}^{2}&=&\sum_{n\in\mathbb{Z}}
w_n\mid f(\mid z_n \mid^2)\mid^2\mid z_n \mid^2 \le
\sum_{n\in\mathbb{Z}}w_n g(\mid z_n\mid^2)^2\mid z_n\mid^2
\nonumber\\
&\le& \sum_{n\in\mathbb{Z}}w_n g(\mid\mid z\mid\mid_{\ell_{w}^{2}}^2)^2\mid z_n\mid^2
\le\{\max_{\rho\in[0,R^2]} g(\rho)\}^2 \sum_{n\in\mathbb{Z}}w_n \mid z_n\mid^2 \le c(R)
\mid\mid z\mid\mid_{\ell_{w}^{2}}. \nonumber
\end{eqnarray}
for some positive constant $c(R)$. Thus we conclude that the operator $T$ is bounded on bounded sets of $\ell_{w}^{2}$.
\\
(ii) Since $\mathrm{(N_2)}$ holds,
for some $\theta \in (0,1)$ we have
\begin{eqnarray}
\mid\mid T(z)-T(z^{'})\mid\mid_{\ell_{w}^{2}} 
&\le& 2\sum_{n\in\mathbb{Z}} w_n  \mid f(\mid z_n\mid^2)\mid^2 \mid z_n -z_n^{'}\mid^{2}\nonumber\\
&&+ 2\sum_{n\in\mathbb{Z}}w_n \mid f^{'} (\theta \mid z_n\mid^2
+(1-\theta)\mid z_n^{'}\mid^2)\mid^{2}(\mid z_n\mid +\mid
z_n^{'}\mid)^{2}  \mid z_n^{'}\mid^{2} \mid z_n -z_n^{'}\mid^2.
\nonumber
\end{eqnarray}
A similar inequality to (\ref{prop4}) should be obtained: We may see that
\begin{eqnarray}
\mid\mid T(z)-T(z^{'})\mid\mid_{l_{w}^{2}}^{2} &\le& 2\sum_{n\in\mathbb{Z}} w_{n}\mid f(\mid\mid  z\mid\mid_{\ell_{w}^{2}}^{2})\mid^{2}\mid z_n -z_{n}^{'}\mid^{2}\nonumber\\
&&+2\{ \max_{\rho \in [0,2R^2]}
g_1(\rho) \}^{2} c(R) \sum_{n\in\mathbb{Z}}w_n \mid
z_{n}-z_{n}^{'}\mid^{2}
\nonumber \\
&\le& c_1(R)( \mid\mid z -z^{'} \mid\mid_{\ell_{w}^{2}}^{2}.
\nonumber
\end{eqnarray}
This concludes the proof of the Lemma.\ \ $\diamond$
\begin{smallTheorem}
\label{thm:LOCALLIPS2}
The operator $A :\ell_{w}^{2}\rightarrow \ell_{w}^{2}$, defined by
$(Au)_n=u_{n+1}-2u_{n}+u_{n-1}$ is globally Lipschitz on $\ell^2_w$.
\end{smallTheorem}
{\bf Proof:}
Let $u,\,v\in B_R$. Then $(Au)_{n\in\mathbb{Z}}-(Av)_{n\in\mathbb{Z}}=(u_{n+1}-v_{n+1})-2(u_n-v_n)+(u_{n-1}-v_{n-1})$ and it follows that $||Au-Av||_{\ell^2_w}\leq 4||u-v||_{\ell^2_w}$.\ $\diamond$\newline

Using the above two Lemmas we may restate the local existence
result, in the case of weighted spaces.
\begin{rigor1}
Let assumption $(\mathrm{N_2})$ be satisfied, assume that $u_0\in\ell^2_w$. There exists
$T^*(u_0)>0$, such that for all $0< T<T^*(u_0)$, there exists a unique solution of the problem (\ref{lat1})-(\ref{lat2}), $u(t)\in \mathrm{C}^1([0,T],\ell^2_w)$.
\end{rigor1}
{\bf Proof:}\ \ This time, we write (\ref{lat1})-(\ref{lat2}), as an ordinary differential equation in $\ell^2_w$
\begin{eqnarray*}
\label{absode}
\dot{u}(t)&+&\Phi(u(t))=0,\\
u(0)&=&u_0,\nonumber
\end{eqnarray*}
where $\Phi(u)=T(u)+A(u)+i\delta u-g$  and $u(t)$ lies in
$\ell^2_w$. Lemmas \ref{thm:LOCALLIPS}-\ref{thm:LOCALLIPS2}
suffice for the application of standard existence and uniqueness
Theorems for ordinary differential equations in Banach spaces
\cite[pg. 78-82]{zei85}. \  $\diamond$
\\ \\
In order to prove the existence of global solutions we need to obtain some {\it a priori} estimates for the solution. In the dissipative case this is achieved by proving the existence of a globally attracting ball of finite radius in $\ell_{w}^{2}$.

\begin{smallTheorem}
\label{ballweighted}
Assume condition $\mathrm{(W)}$ on the weight function and that the damping coefficient satisfies
\begin{eqnarray}
\label{condidamp}
\frac{\delta}{2}-2d_1\underline{d}_2^{-1/2} \ge 0.
\end{eqnarray}
Let $(u_{0,n})_{n\in\mathbb{Z}}=u_0, (g_n)_{n\in\mathbb{Z}}=g\in\ell^2_w$ and condition $\mathrm{(N_2)}$ be satisfied.  A dynamical system can be defined by
(\ref{lat1})-(\ref{lat2}),
\begin{eqnarray}
\label{dynamicalw} S(t):u_0\in {\ell}^2_w\rightarrow u(t)\in
{\ell}^2_w,
\end{eqnarray}
possessing a bounded absorbing set $\mathcal{B}_0$ in ${\ell}^2_w$:
For every bounded set $\mathcal{B}$ of ${\ell}^2_w$, there exists
$t_0(\mathcal{B},\mathcal{B}_0)$ such that for all  $t\geq
t_0(\mathcal{B},\mathcal{B}_0)$, it holds
$S(t)\mathcal{B}\subset\mathcal{B}_0$.
\end{smallTheorem}
{\bf Proof:} We multiply (\ref{lat1}) with $w_n\overline{u}_n$,  $n\in\mathbb{Z}$ add over all
lattice sites and keep the imaginary part. Working similarly as in
Lemma \ref{DNLStail}, we obtain
\begin{eqnarray}
\frac{1}{2}\frac{d}{dt}\sum_{n\in\mathbb{Z}}w_n \mid u_n\mid^2
&-&\sum_{n\in\mathbb{Z}}(w_{n+1}-w_{n})(u_{1,n}u_{2,n+1}-u_{2,n}u_{1,n+1})+
\delta \sum_{n\in\mathbb{Z}}w_{n}\mid u_{n}\mid^2 \nonumber\\
&\le& \frac{\delta}{2}\sum_{n\in\mathbb{Z}}w_n \mid u_n\mid^2
+\frac{1}{2\delta}\sum_{n\in\mathbb{Z}}w_n \mid g_n\mid^2
\nonumber
\end{eqnarray}
Using the assumptions on the weight function $\mathrm{(W})$, we find that
\begin{eqnarray}
\sum_{n\in\mathbb{Z}}\mid(w_{n+1}-w_{n})(u_{1,n}u_{2,n+1}-u_{2,n}u_{1,n+1})\mid
&\leq&  d_1 \sum_{n\in\mathbb{Z}}w_n
\mid(u_{1,n}u_{2,n+1}-u_{2,n}u_{1,n+1})\mid
\nonumber \\
&\leq& d_1 \left(\sum_{n\in\mathbb{Z}} w_{n} \mid u_{1,n} u_{2,n+1}\mid
+  \sum_{n\in\mathbb{Z}} w_{n}\mid u_{2,n} u_{1,n+1}\mid\right)
\end{eqnarray}
We now estimate the sums in the above inequality as follows:
\begin{eqnarray}
 \sum_{n\in\mathbb{Z}} w_n \mid u_{1,n}u_{2,n+1} \mid &=& 
\sum_{n\in\mathbb{Z}} w_n^{1/2} w_n^{1/2} \mid u_{1,n}u_{2,n+1} \mid
\le  \left(\sum_{n\in\mathbb{Z}} w_{n}\mid u_{1,n}\mid^2\right)^{1/2} 
\left(\sum_{n\in\mathbb{Z}} w_n \mid u_{2,n+1}\mid^2\right)^{1/2}
\nonumber \\
&\le& \underline{d}_2^{-1/2} 
\left(\sum_{n\in\mathbb{Z}}w_{n}\mid u_{1,n}\mid^2\right)^{1/2} 
\left(\sum_{n\in\mathbb{Z}}w_{n+1}\mid u_{2,n+1}\mid^2\right)^{1/2} \le
\underline{d}_2^{-1/2} \mid\mid u \mid\mid^2_{\ell_{w}^{2}},
\nonumber
\end{eqnarray}
and similarly,
\begin{eqnarray}
\sum_{n\in\mathbb{Z}} w_n \mid u_{2,n}u_{1,n+1} \mid \le
\underline{d}_2^{-1/2} \mid\mid u \mid\mid_{\ell_{w}^{2}}^2. \nonumber
\end{eqnarray}
Therefore,
\begin{eqnarray}
\sum_{n\in\mathbb{Z}}\mid (w_{n+1}-w_{n})(u_{1,n}u_{2,n+1}-u_{2,n}u_{1,n+1})\mid
\leq 2 d_1 \underline{d}_2^{-1/2} \mid\mid u
\mid\mid_{\ell_{w}^{2}}^2. \nonumber
\end{eqnarray}
Thus, we obtain the differential inequality
\begin{eqnarray}
\frac{1}{2}\frac{d}{dt}\mid\mid u \mid\mid_{\ell_{w}^{2}}^{2}+\left(\frac{\delta}{2}-2d_1\underline{d}_{2}^{-1/2}
\right) \mid\mid u\mid\mid_{\ell_{w}^{2}}^{2} \le \frac{1}{2\delta} \mid\mid g\mid\mid_{\ell_{w}^{2}}^{2}
\nonumber
\end{eqnarray}
Using the Gronwall inequality, if
$\frac{\delta}{2}-2d_1\underline{d}_2^{-1/2} \ge 0$ we obtain the
existence of the attracting ball in $\ell_{w}^{2}$. This concludes
the proof of the Proposition. $\diamond$
\begin{remark} {\em In the case of exponential weight $w_{n}=exp(\lambda n)$, it appears that $d_1=exp(\lambda)-1$ and $\underline{d}_2=exp(\lambda)$. The condition for existence of attracting ball  becomes $8 sinh(\lambda /2) \le \delta$.}
\end{remark}

We may further prove the existence of a global attractor
attracting all bounded sets of $\ell^2_w$. The next Lemma provides us
with tail estimates in the weighted space.

\begin{smallTheorem}
Let $(u_{0,n})_{n\in\mathbb{Z}}=u_0\in\mathcal{B}$ where $\mathcal{B}$ is a bounded set of $\ell^2_w$ and $(g_n)_{n\in\mathbb{Z}}=g\in\ell^2_w$. Under assumptions $\mathrm{(N_2)}$ on the nonlinearity, $\mathrm{(W)}$ on the weight function and  assumption (\ref{condidamp}) on the dissipation, for any  $\eta >0$, there exist $T(\eta )$ and $K(\eta )$ such that the solution
$(u_{n})_{n\in\mathbb{Z}}=u$ of (\ref{lat1})-(\ref{lat2}) satisfies  for all $t\geq T(\eta )$, the estimate
\begin{eqnarray}
\label{precomw}
\sum_{\mid n\mid >2M}w_n\mid u_n\mid^2 \le \frac{2\eta}{\delta},\;\;\mbox{for any}\;\;
M>K(\eta ).
\end{eqnarray}
\end{smallTheorem}
{\bf Proof:} The proof follows closely the proof of Lemma
\ref{DNLStail} for the corresponding tail estimates in the space
$\ell^{2}$, only that we now multiply the equation by $w_n \theta_n \overline{u}_n$, instead of multiplying simply by $\theta_n
\overline{u}_n$. 
We obtain the following inequality
\begin{eqnarray}
\label{te1}
\frac{1}{2}\frac{d}{dt}\sum_{n\in\mathbb{Z}}w_n\theta_n \mid u_n\mid^2
-I_2+
\delta \sum_{n\in\mathbb{Z}}w_{n}\theta_n\mid u_{n}\mid^2
\le \frac{\delta}{2}\sum_{n\in\mathbb{Z}}w_n\theta_n \mid u_n\mid^2
+\frac{1}{2\delta}\sum_{n\in\mathbb{Z}}w_n\theta_n \mid g_n\mid^2,
\end{eqnarray}
where $$I_2=\sum_{n\in\mathbb{Z}}\{(\theta_{n+1}w_{n+1}-\theta_n w_n)(u_{1,n}u_{2,n+1}-u_{2,n}u_{1,n+1})\}.$$
We rewrite
$$ (\theta_{n+1}w_{n+1}-\theta_n w_n)=(\theta_{n+1}-\theta_n)w_{n+1}+\theta_n (w_{n+1}-w_{n}).$$  
Using arguments similar to those used in the proof of Lemma \ref{DNLStail} and Lemma \ref{ballweighted} we see that
\begin{eqnarray}
\label{teA}
\sum_{n\in\mathbb{Z}}|\theta_{n+1}-\theta_n|\,|w_{n+1}|\,|u_{1,n}u_{2,n+1}-u_{2,n}u_{1,n+1}|\leq \frac{C_2\rho_2^2}{M},
\end{eqnarray}
where $C_2$ depends only on $d_1,\underline{d}_2$ and $\rho_2$ denotes the absorbing ball in $\ell^2_w$.
We also have the inequality
\begin{eqnarray}
\label{teB1}
\sum_{n\in\mathbb{Z}}\theta_n\,|w_{n+1}-w_n|\,|u_{1,n}u_{2,n+1}-u_{2,n}u_{1,n+1}|\leq
d_1\sum_{n\in\mathbb{Z}}\theta_nw_n|u_{1,n}u_{2,n+1}|+d_1\sum_{n\in\mathbb{Z}}\theta_nw_n|u_{2,n}u_{1,n+1}|.
\end{eqnarray}
For the first term of the rhs of (\ref{teB1}), holds the estimate
\begin{eqnarray}
\label{teB2}
d_1\sum_{n\in\mathbb{Z}}\theta_nw_n|u_{1,n}u_{2,n+1}|&\leq& d_1\left(\sum_{n\in\mathbb{Z}}\theta_nw_n|u_{1,n}|^2\right)^{1/2}
\left(\sum_{n\in\mathbb{Z}}\theta_nw_n|u_{2,n+1}|^2\right)^{1/2}\nonumber\\
&\leq& d_1\underline{d}_2^{-1/2}
\left(\sum_{n\in\mathbb{Z}}\theta_nw_n|u_{1,n}|^2\right)^{1/2}
\left(\sum_{n\in\mathbb{Z}}\theta_nw_{n+1}|u_{2,n+1}|^2\right)^{1/2}.
\end{eqnarray}
An application of Young's inequality to the rhs of (\ref{teB2}), implies that
\begin{eqnarray}
\label{teB3}
d_1\sum_{n\in\mathbb{Z}}\theta_nw_n|u_{1,n}u_{2,n+1}|&\leq&
\frac{1}{2}d_1\underline{d}_2^{-1/2}\sum_{n\in\mathbb{Z}}\theta_nw_n|u_{1,n}|^2
+\frac{d_1^2\underline{d}_2^{-1}}{2d_1\underline{d}_2^{-1/2}}
\sum_{n\in\mathbb{Z}}\theta_nw_{n+1}|u_{2,n+1}|^2\nonumber\\
&\leq&\frac{1}{2}d_1\underline{d}_2^{-1/2}\sum_{n\in\mathbb{Z}}\theta_nw_n|u_{1,n}|^2
+\frac{1}{2}d_1\underline{d}_2^{-1/2}\sum_{n\in\mathbb{Z}}\theta_{n+1}w_{n+1}|u_{2,n+1}|^2\nonumber\\
&&+\frac{1}{2}d_1\underline{d}_2^{-1/2}\sum_{n\in\mathbb{Z}}(\theta_n-\theta_{n+1})w_{n+1}|u_{2,n+1}|^2\nonumber\\
&\leq&
d_1\underline{d}_2^{-1/2}\sum_{n\in\mathbb{Z}}\theta_nw_n|u_{1,n}|^2+\frac{C_3\rho_2^2}{M},
\end{eqnarray}
and similarly for the second term of the rhs of (\ref{teB1})
\begin{eqnarray}
\label{teB4}
d_1\sum_{n\in\mathbb{Z}}\theta_nw_n|u_{2,n}u_{1,n+1}|\leq d_1\underline{d}_2^{-1/2}\sum_{n\in\mathbb{Z}}\theta_nw_n|u_{2,n}|^2+\frac{C_4\rho_2^2}{M}.
\end{eqnarray}
Using (\ref{te1})-(\ref{teB4}) we derive that
\begin{eqnarray*}
\frac{1}{2}\frac{d}{dt}\sum_{n\in\mathbb{Z}}w_n\theta_n \mid u_n\mid^2+
\left(\frac{\delta}{2}-2d_1\underline{d}_2^{-1/2}\right)\sum_{n\in\mathbb{Z}}w_n\theta_n \mid u_n\mid^2
\le +\frac{C\rho_2^2}{M}
+\frac{1}{2\delta}\sum_{n\in\mathbb{Z}}w_n\theta_n \mid g_n\mid^2.
\end{eqnarray*}
with $C$ depending only on $d_1,\underline{d}_2$. The rest follows by direct generalization of the arguments of
Lemma \ref{DNLStail}. $\diamond$
\\ \\
In complete analogy as before we may prove the asymptotic compactness of the semigroup in the weighted spaces,
\begin{rigor2}
\label{asymcompweighted}
The semigroup $S(t)$ is asymptotically compact in $\ell^{2}_w$, that is if the sequence $\phi_n$ is bounded in $\ell^{2}_w$ and $t_n \rightarrow \infty$ then $S(t_n)\phi_n$ is precompact in $\ell^2_w$.
\end{rigor2}
This result leads to the existence of a global attractor in the weighted spaces.
\begin{rigor1}
\label{dynamics1w} Let conditions $\mathrm{(N_2),(W)}$ and
(\ref{condidamp}) be fulfilled. The
semigroup $S(t)$ associated to (\ref{lat1})-(\ref{lat2}) possesses
a global attractor
$\mathcal{A}=\omega(\mathcal{B}_0)\subset\mathcal{B}_0\subset\ell^2_w$
which is compact, connected and maximal among the functional
invariant sets in $\ell^2_w$.
\end{rigor1}
\section{Remarks on possible extensions}
We conclude, by mentioning some other examples of DNLS-type equations  for which, extensions of the results on the global solvability and the existence of global attractors, could be investigated. 

A first  example  is provided by  the DNLS with potential \cite{ANY1,ANY2},  i.e (\ref{lat1})-(\ref{lat2}) with a nonlinearity of the form $F(u)_{n\in\mathbb{Z}}=F(u_{n})+V(n)u$, where
$V:\mathbb{R}\rightarrow\mathbb{R}$ is a real valued potential $|V(s)|\leq c$ for every $s\in\mathbb{R}$. The potential function expresses the inhomogeneity properties of the medium.

Another interesting DNLS model, may be given by the spatial discretization of a modified NLS equation with viscocity
\begin{eqnarray}
\label{pdeM}
iu_t+\frac{1}{2}u_{xx}-|u|^2u=i\delta u_{xx},\;\;\delta >0
\end{eqnarray}
which  describes the light propagation in an array of optical fibers in a weakly lossy medium. Taking into account some external excitation and considering a general nonlinearity, the discretized counterpart of (\ref{pdeM}) reads as (see \cite{Mal1})
\begin{eqnarray}
\label{dpdeM}
i\dot{u}_n+(1-i\delta)(u_{n+1}-2u_n+u_{n-1})+f(|u_n|^2)u_n=g_n.
\end{eqnarray}

Typical DNLS-models include equations of the form (\cite{Gupta,Kim})
\begin{eqnarray*}
\label{dpdeM2}
i\dot{U}_n+(\alpha_{n+1}U_{n+1}+\beta_{n+1}U_{n+1})+i\delta_nU_n+\gamma_nf(|U_n|^2)U_n=g_n.
\end{eqnarray*}
For the coupling and nonlinear strength and dissipation parameter $\alpha_n,\beta_n,\gamma_n,\delta_n\in\mathbb{R}$ we assume that there exist some  constants $c_1,c_2, c_3$ such that
\begin{eqnarray*}
\label{strc}
|\alpha_n|,|\beta_n|,|\gamma_n|\leq c_1,\;\;c_2\leq\delta_n\leq c_3,\;n\in\mathbb{Z}.
\end{eqnarray*}

The validity of the results could be examined for DNLS-type equations  considered in $\mathbb{Z}^N$, $N>1$,
\begin{eqnarray}
\label{extend1}
i\dot{u}_n+(Au)_n &+&i\delta_n u_n+f(|u_n|^2)u_n=g_n,\\
\label{lat2m}
u_n(0)&=&u_{n,0},\;\;n=(n_1,n_2,...n_N)\in\mathbb{Z}^N.
\end{eqnarray}
where dissipation satisfies (\ref{strc}) and the linear operator $A$ has a decomposition introduced in (\cite{SZ2}): It is assumed that
\begin{eqnarray}
\label{gop1}
A=A_1+A_2+...+A_N
\end{eqnarray}
and
\begin{eqnarray}
\label{gop2}
A_j=B_j^*B_j=B_j^*B_j,\;\;||B_j||_{\mathcal{L}(\ell^2)}\leq M,\;\;j=1,2,...,N,
\end{eqnarray}
for some bounded linear operators $B_j:\ell^2\rightarrow\ell^2$,
and its adjoint $B_j^*$, defined as
\begin{eqnarray}
\label{gop3}
(B_ju)_n&=&\sum_{l=-m_0}^{l=m_0}C_{j,l}u_{n_{jl}},\;\;\mbox{for all}\;\;u=(u_n)_{n\in\mathbb{Z}^N},\;\;j=1,2,...N\\
(B_j^*u)_n&=&\sum_{l=-m_0}^{l=m_0}C_{j,-l}u_{n_{jl}},\;\;n_{jl}=(n_1,n_2,...,n_{j-1},n_j+l,n_{j+1},...,n_N)\in\mathbb{Z}^N.
\end{eqnarray}
Under assumptions (\ref{gop1})-(\ref{gop3}) similar relations to (\ref{diffop2})-(\ref{diffop3}) are satisfied,
\begin{eqnarray*}
(Bu,v)_{\ell^2}=(u,B^*v)_{\ell^2},\;\;
(Au,v)_{\ell^2}=-(Bu,Bv)_{\ell^2},\;\;u,v\in\ell^2.
\end{eqnarray*}
Especially for the conservative case, $\delta_n=0$, $n\in\mathbb{Z}^N$, and the power-law nonlinearity, (\ref{gop1})-(\ref{gop3}) are sufficient for the extension of Theorem \ref{global} to (\ref{extend1})-(\ref{lat2m}).
Although this is a  simple observation, it is important,  since it is a rigorous verification that spatial dimension does not play a particular role on the global existence of the solutions as it happens in the continuous model \cite{bang}, \cite{cazh,cazS}. For a discussion on the existence of global attractors for (\ref{extend1})-(\ref{lat2m}), one could be based on the framework provided by \cite{SZ2}.
\\ \\
{\bf Acknowledgements} We would like also  to thank Professors D.  Frantzeskakis, K. Hizanidis, and N. Stavrakakis for helpful discussions. 
%\bibitem[ ]

\end{document}